\documentclass[a4paper,10pt]{amsart}
\usepackage{amssymb,amsmath,amsfonts,amsthm}
\usepackage[dvips]{graphicx}
\usepackage{psfrag}
\numberwithin{equation}{section}\voffset 0cm \hoffset -1cm
\newtheorem{thm}{Theorem}[section]
\newtheorem{lem}[thm]{Lemma}
\newtheorem{pro}[thm]{Proposition}
\newtheorem{dfn}[thm]{Definition}

\newtheorem{cor}[thm]{Corollary}
\newtheorem{rk}[thm]{Remark}

\def\bb{\begin}               \def\cal{\mathcal}
           \def\ea{\end{array}}
          \def\ec{\end{center}}
     \def\ed{\end{description}}
        \def\ee{\end{equation}}
       \def\eea{\end{eqnarray}}
     \def\eeaa{\end{eqnarray*}}
\def\bt{\bb{thebibliography}} \def\et{\end{thebibliography}}
\def\bib{\bibitem}
\def\aaa{\vdash\hskip-2pt\dashv}
\def\h{\pitchfork}

\def\q{\quad}

\def\Pf {\noindent {\bf Proof}\q}
\def\qed {\hfill $\Box$\vskip5pt}

\title{Newhouse phenomenon and homoclinic classes}
\author{Jiagang Yang}
\date{\today}
\thanks{Partially supported by TWAS-CNPq, FAPERJ}

\address{IMPA, Est. D. Castorina 110, 22460-320 Rio de Janeiro, Brazil}
\email{yangjg\@@impa.br}

\begin{document}

\begin{abstract}
We show that for a $C^1$ residual subset of diffeomorphisms far
away from tangency, every non-trivial chain recurrent class that
is accumulated by sources ia a homoclinic class contains periodic
points with index 1 and it's the Hausdorff limit of a family of
sources.
\end{abstract}
\maketitle \tableofcontents
\section{Introduction}
In the middle of last century, with many remarkable work,
hyperbolic diffeomorphisms have been understood very well, but
soon people discovered that the set of hyperbolic
 diffeomorphisms are not dense among differential dynamics, two kinds of counter
 examples were described, one associated with heterdimension cycle was given by
 R.Abraham and Smale \cite{1AS} and then given by Shub \cite{1Sh}  and Ma\~{n}\'{e} \cite{1M2}, another counter
 example associated with homoclinic tangency was given by Newhouse \cite{1N1} \cite{1N2}. In fact,
 Newhouse got an open set $\cal U$$\subset C^2(M)$ where $dim(M)=2$ such that there exists a
 $C^2$ generic subset $R\subset \cal U$ and for any $f\in R$ , $f$ has infinite
 sinks or sources. Such complicated phenomena (there exist an open set $\cal U$ in $C^r(M)$
 and a generic subset $R\subset \cal U$, such that any $f\in R$ has infinite sinks or sources)
 is called $C^r$ Newhouse phenomena today, and we say $C^r$ Newhouse phenomena happens at $\cal U$.

In last 90's, some new examples of Newhouse phenomena were found,
\cite{1PV} generalized Newhouse phenomena
 to high dimensional manifold ($dimM>2$) but with the same topology $C^r(r>1)$. \cite{1BD} used a new tool 'Blender'
 to show the existence of $C^1$ Newhouse phenomena on manifold with $dim(M)>2$. Until now, all the construction
 of $C^r$ Newhouse phenomena relate closely with homoclinic tangency, more precisely, all the open set
 $\cal U$
given by the construction above which happens Newhouse phenomena
there will have ${\cal U}\subset \overline{HT}$. We hope that it's a
necessary condition for $C^r$ Newhouse phenomena happens at $\cal
U$. Pujals states it
as a conjecture.\\

\noindent{\bf Conjecture\ }(Pujals): {\it If $C^r$ Newhouse
phenomena happens at $\cal U$, then $\cal U$ is contained in
$\overline{HT^r}$.}\\

  When $r=1$ and $M$ is a compact surface, with Ma\~{n}\'{e}'s work \cite{1M3},
  Pujals' conjecture is equivalent with the
famous $C^1$ Palis strong conjecture. \\

\noindent{\bf $C^1$ Palis strong conjecture\ }: {\it
Diffeomorphisms of $M$ exhibiting either a homoclinic tangency or
heterodimensional
cycle are $C^1$ dense in the complement of the $C^1$ closure of hyperbolic systems.}\\

  In the remarkable paper \cite{1PS1} they proved $C^1$ Palis strong conjecture on $C^1(M)$ when $M$ is a boundless compact
surface, so in such case Pujals' conjecture is right. In
\cite{1PS2} they gave many relations between $C^2$ Newhouse
phenomena and $\overline{HT^1}$. In this paper we just consider
$C^1$ Newhouse phenomena, and we show that if $C^1$ Newhouse
phenomena happens in an open set $\cal U$ $\subset
C^1(M)\backslash \overline{HT^1}$, it should have some special
properties, in fact, in \cite{1BD} they found an open set $\cal U$
$\subset \overline {(HT^1)}$ and there exists a generic subset
$R\subset \cal U$ such that any $f\in R$ has infinite sinks or
sources stay near a chain recurrent class, and such class does not
contain any periodic points, such kind of chain recurrent class is
called aperiodic class now. Here we proved that in
$\overline{HT}^c$, if there exists Newhouse phenomena, the sinks
or sources will just stay
near a special kind of homoclinic class.\\

\noindent{\bf Theorem 1} {\it There exists a generic subset $R
\subset C^1(M)\backslash \overline{HT^1}$, such that for $f\in R$
and $C$ is any non-trivial chain recurrent class of $f$, if
$C\bigcap P^*_0\neq \phi$, $C$ should be a homoclinic class
containing index 1
periodic points and $C$ is an index 0 fundamental limit.}\\

 Theorem 1 means that if we want to disprove the existence of
Newhouse phenomena in $C^1(M)\setminus\overline{HT}$, we just need
study the homoclinic class containing index 1 periodic point.

In $\S$3 we'll state some generic properties. In
$\S$4 we'll introduce a special minimal non-hyperbolic set and
theorem 1 will be proved in $\S$5.

\noindent{\bf Acknowledgements:} This paper is part of the
author's thesis, I would like to thank my advisor Professor
Marcelo Viana for his support and enormous encouragements during
the preparation of this work. I would like to thank Professor
Shaobo Gan for checking the details of the proof and finding out
an essential gap in the original argument which is crucial to the
work. I also thank Professors Jacob Palis, Lan Wen, Enrique R.
Pujals, Lorenzo Diaz, Christian Bonatti for very helpful remarks.
Finally I wish to thank my wife, Wenyan Zhong, for her help and
encouragement.

\section{Notations and definitions}

  Let $M$ be a compact boundless Riemannian manifold, since when
  $M$ is a surface \cite{1PS1} has proved that hyperbolic diffeomorphisms
  are open and dense in $C^1(M)\setminus \overline{HT}$, we suppose $dim(M)=d>2$ in this paper.
  Let $Per(f)$ denote the set of periodic points of $f$ and $\Omega (f)$ the
  non-wondering set of $f$, for $p \in Per(f)$, $\pi (p)$ means the period of $p$. If
  $p$ is a hyperbolic periodic point, the index of $p$ is the dimension
  of the stable bundle. We denote $Per_i(f)$ the set of the index $i$ periodic points of $f$, and we call a point $x$ is an index $i$ preperiodic
  point of $f$ if there exists a family of diffeomorphisms $g_n\stackrel{C^1}{\longrightarrow}f$, where $g_n$ has an index $i$ periodic point $p_n$ and ${p_n}\longrightarrow x$.
  $P^*_i(f)$ is the set of index $i$ preperiodic points of $f$, it's easy to know $\overline{P_i(f)}\subset P^*_i(f)$.

  Let $\Lambda$ be an invariant compact set of $f$, we say $\Lambda$ is an index $i$ fundamental
  limit if there exists a family of diffeomorphisms ${g_n}$ $C^1$ converging to $f$,
  $p_n$ is an index $i$ periodic point of $g_n$ and $Orb(p_n)$ converge to $\Lambda$ in \emph{Hausdorff}
  topology. So if $\Lambda (f)$ is an index $i$ fundamental limit, we have $\Lambda (f)\subset P^*_i(f)$.

  For two points $x,y\in M$ and some $\delta >0$, we say there exists a $\delta$-pseudo orbit
  connects $x$ and $y$ means that there exist points $x=x_0,x_1,\cdots,x_n=y$ such
  that $d(f(x_i),x_{i+1})<\delta$ for $i=0,1,\cdots ,n-1$, we denote it $x \underset{\delta}{\dashv}y$. We say $x\dashv y$ if
  for any $\delta >0$ we have $x \underset{\delta}{\dashv} y$
  and denote $x \aaa y$ if $x\dashv y$ and $y\dashv x$. A point $x$ is called a chain recurrent point if $x\aaa x$. $CR(f)$ denotes
  the set of chain recurrent points of $f$, it's easy to know that $\aaa$ is an
  closed equivalent relation on $CR(f)$, and every equivalent class of such
  relation should be compact and is called chain recurrent class. Let $K$ be
  a compact invariant set of $f$, if $x,y$ are two points in $K$, we'll denote $x \underset{K}{\dashv} y$
  if for any $\delta >0$, we have a $\delta$ -pseudo orbit in $K$ connects $x$ and $y$. If for any
  two points $x,y\in K$ we have $x \underset{K}{\dashv} y$, we call $K$ a chain recurrent set. Let
  $C$ be a chain recurrent class of $f$, we call $C$ is an aperiodic
class if $C$ does not contain periodic point.

  Let $\Lambda$ be an invariant compact set of $f$, for $0<\lambda <1$ and $1\leq i<d$, we say $\Lambda$ has
  an index $i-(l,\lambda)$ dominated splitting if we have a continuous invariant
  splitting $T_{\Lambda}M=E\oplus F$ where $dim(E_x)=i$ for any $x\in \Lambda$ and $\parallel Df^l|_E(x)\parallel\cdot
  \parallel Df^{-l}|_F(f^lx)\parallel <\lambda$ for all $x\in \Lambda$.
  For simplicity, sometimes we just call $\Lambda (f)$ has an index $i$ dominated splitting.
  A compact invariant set can have many dominated splittings, but for fixed $i$,
  the index $i$ dominated splitting is unique.

  We say a diffeomorphism $f$ has $C^r$ tangency if $f\in C^r(M)$, $f$ has hyperbolic
  periodic point $p$ and there exists a non-transverse intersection between
  $W^s(p)$ and $W^u(p)$. $HT^r$ is the set of the diffeomorphisms which have
  $C^r$ tangency, usually we just use $HT$ denote $HT^1$. We call a
diffeomorphism $f$
  is far away from tangency if $f\in C^1(M)\setminus \overline{HT}$. The following proposition shows the
  relation between dominated splitting and far away from
  tangency.

\begin{pro}\label{2.1} (\cite{1W1}) $f$ is $C^1$ far away from tangency if
and only if there exists $(l,\lambda)$ such that $P^*_i(f)$ has
index $i-(l,\lambda)$ dominated splitting for $0<i<d$.\end{pro}

  Usually dominated splitting is not a hyperbolic splitting, Ma\~{n}\'{e} showed
  that in some special case, one bundle of the dominated splitting is
  hyperbolic.

\begin{pro}\label{2.2}(\cite{1M3}) Suppose $\Lambda (f)$ has an index $i$
dominated splitting $E\oplus F\  (i\neq 0)$, if $\Lambda
(f)\bigcap P^*_j(f)=\phi$ for $0\leq j<i$, then $E$ is a
contracting bundle.\end{pro}

\section{Generic properties}
 For a topology space $X$, we call a set $R\subset X$ is a generic subset
  of $X$ if $R$ is countable intersection of open and dense subsets of $X$,
  and we call a property is a generic property of $X$ if there exists
  some generic subset $R$ of $X$ holds such property. Especially, when
  $X=C^1(M)$ and $R$ is a generic subset of $C^1(M)$, we just call $R$ is
  $C^1$
  generic, and we call any generic property of $C^1(M)$ 'a $C^1$ generic
  property' or 'the property is $C^1$ generic'.

  Here we'll state some well known $C^1$ generic properties.

\begin{pro}\label{3.1} There is a $C^1$ generic subset $R_0$ such
that for any $f\in R_0$, one has
\begin{itemize}
  \item[1)]  $f$ is Kupka-Smale (every periodic
point $p$ in $Per(f)$ is hyperbolic and the invariant manifolds of
periodic points are everywhere transverse).
  \item[2)] $CR(f)=\Omega =\overline{Per(f)}$.
  \item[3)] $P^*_i(f)=\overline{P_i(f)}$
  \item[4)]  any chain recurrent set is the \textit{Hausdorff}
limit of periodic orbits.
  \item[5)]  any index $i$ fundamental limit is the
Hausdorff limit of index $i$ periodic orbits of $f$.
  \item[6)] any chain
recurrent class containing a periodic point $p$ is the homoclinic
class $H(p,f)$.
  \item[7)] Suppose $C$ is a homoclinic class of $f$, and $j_0=min\{j:C\bigcap Per_j(f)\neq
  \phi\},\ j_1=max\{j: C\bigcap Per_j(f)\neq \phi\},$
then for any $j_0\leq j\leq j_1$, we have $C\bigcap Per_j(f)\neq
\phi$.
\end{itemize}
\end{pro}
  By proposition \ref{3.1}, for any $f$ in $R_0$, every chain recurrent class
  $C$
  of $f$ is either an aperiodic class or a homoclinic class. If $\# C=\infty$,
  we call $C$ is non-trivial.

Let $R=R_0\setminus \overline{HT}$, we'll show that the generic subset $R$ of $\overline{HT}^c$ will satisfy theorem 1.

\section{A special minimal set}
  Let $f\in R$, $C$ is a non-trivial chain recurrent class of $f$, and
  $j_0=\min\{j:C\bigcap P^*_j\neq \phi\}$.

\begin{dfn}\label{4.1}: An invariant compact subset $\Lambda$ of $f$ is
called minimal if all the invariant compact subsets of $\Lambda$
are just $\Lambda$ and $\phi$. An invariant compact subset
$\Lambda$ of $f$ is called minimal index $j$ fundamental limit if
$\Lambda$ is an index $j$ fundamental limit and any invariant
compact subset $\Lambda_0\varsubsetneq \Lambda$ is not an index
$j$ fundamental limit.\end{dfn}

\begin{lem}\label{4.2} If $C\bigcap P^*_j\neq \phi$, there always exists a minimal index $j$ fundamental limit in
$C$.\end{lem}

\Pf Let $H=\{\tilde{\Lambda}: \tilde{\Lambda}\subset C$ is an index $j$
fundamental limit$\}$ and we order $H$ by inclusion. Suppose $x\in
C\bigcap P^*_j$, then there exist $g_n\overset{C^1}{\longrightarrow} f$, $p_n$ is
index $j$ periodic point of $g_n$ and $p_n\longrightarrow x$.
Denote $\Lambda_x=\lim{Orb(P_n)}$, then $\Lambda_x$ is an index
$j$ fundamental limit. It's easy to know $\Lambda_x$ is a chain
recurrent set and $\Lambda_x\subset C$, so $\Lambda_x\in H$. It
means $H\neq \phi$.

  Let $H_\Gamma=\{\Lambda_\lambda:\lambda \in \Gamma\}$ be a totally ordered chain of $H$. Then
  $\Lambda_\infty=\bigcap_{\lambda \in \Gamma}\Lambda_\lambda$ is a compact invariant
   set, in fact, there exists $\{\lambda_i\}^\infty_{i=1}$ such that
   $\Lambda_{\lambda_i}\supset \Lambda_{\lambda_{i+1}}$ and
   $\Lambda_\infty=\bigcap^\infty_{i=1}\Lambda_{\lambda_i}$.

  We claim that $\Lambda_\infty$ is an index $j$ fundamental limit also.\\

\noindent{\bf Proof of the claim\q} From generic property $5)$ of
proposition \ref{3.1} and $f\in R$, for any $\varepsilon
>0$, there exists periodic point $p_i$ such that $p_i\in Per_j(f)$ and
$d_H(Orb(p_i),\Lambda_{\lambda_i})<\frac{\varepsilon}{2}$. When
$i$ is big enough, we'll have
$d_H(\Lambda_{\lambda_i},\Lambda_\infty)<\frac{\varepsilon}{2}$,
so for any $\varepsilon >0$, there exists $p_i\in Per_j(f)$such
that $d_H(Orb(p_i),\Lambda_\infty)< \varepsilon$. \qed

  Now by Zorn's lemma, there exists a minimal index $j$ fundamental limit in $C$.        \qed

  Suppose $\Lambda$ is a minimal index $j_0$ fundamental limit of $C$, the main aim of this section is the following lemma.

\begin{lem} \label{4.3} Suppose $f\in R$, $C$ is a non-trivial chain recurrent class of
$f$, $j_0=\min\{j:C\bigcap P^*_j\neq \phi\}$. Let $\Lambda$ be any
minimal index $j_0$ fundamental limit in $C$, then
\begin{itemize}
  \item[a)]
       either $\Lambda$ is a non-trivial minimal set with partial hyperbolic splitting
       $T|_\Lambda M=E^s_{j_0}\oplus E^c_1\oplus E^u_{j_0+2}$,
  \item[b)]     or $C$ contains a periodic point with index $j_0$ or $j_0+1$ and $C$ is an index $j_0$ fundamental limit.\\
\end{itemize}
\end{lem}
  We postpone the proof of lemma \ref{4.3} to $\S$4.4, before that, I'll give or introduce some results at first.
   In $\S$4.1 I'll give a proof of Shaobo Gan's lemma, in $\S$4.2 I'll introduce Liao's selecting
   lemma and prove a weakly selecting lemma, in $\S$4.3 I'll introduce a powerful tool 'transition' given by [BDP].

\subsection{Shaobo Gan's lemma}

  Let $GL(d)$ be the group of linear isomorphisms of $R^d$, we call
$\xi$ a periodic
  sequence of linear map if $\xi:Z\longrightarrow Gl(d)$ is a sequence of isomorphisms of $R^d$ and there
   exists $n_0\geq 1$ such that $\xi_{j+n_0}=\xi_j$ for all $j$. We denote
   $\pi(\xi)=min\{ n: \xi_{j+n}=\xi_j$ for all $j\}$ the period of $\xi$,
   and we call $\xi$ has index $i$ if the map $\prod \limits^{\pi(\xi)-1}_{j=0}\xi_j$ is hyperbolic and has index $i$,
    we say $\xi$ is contracting if $\xi$ has index $d$. We denote $E^{s(u)}$ the stable
    (unstable) bundle of $\xi$.

  Suppose $\eta$ is a periodic sequence of linear maps also, we call $\eta$ is an $\varepsilon$-perturbation
  of $\xi$ if $\pi(\eta)=\pi(\xi)$ and $\parallel \eta_j-\xi_j\parallel\leq \varepsilon$ for any $j$.

  Let $\{\xi^{\alpha}\}_{\alpha \in \cal A}$ be a family of periodic sequence of linear maps with index $i$, we call
  it is bounded if there exists $K>0$ such that for any $\alpha\in \cal A$ and any $j\in \mathbb{Z}$,
  we have $\parallel \xi^{(\alpha)}_j\parallel<K$. For a family of bounded periodic sequences of linear
  maps $\{\xi^{\alpha}\}_{\alpha \in \cal A}$, we say it's index stable if $\xi^{(\alpha)}$ has index $i$ for all $\alpha\in \cal A$, and there
  exists $\varepsilon_0>0$ such that $\# \{\alpha|$ there exists
  $\eta^{(\alpha)}$ is $\varepsilon_0$-perturbation of $\xi^{(\alpha)}$ and $\eta^{(\alpha)}$ has index
  different with $i\}<\infty$.
  Especially, if it's index $d$ stable, we call $\xi^{(\alpha)}|_{\alpha \in \cal A}$ is uniformly contracting.

  Suppose $f\in C^1(M)$ and $\{p_n(f)\}$ is a family of hyperbolic periodic
  points of $f$ with index $i$, we say ${p_n(f)}$ is index $i$ stable if $\{Df|_{Orb(p_n)}\}^\infty_{n=1}$
  is index $i$ stable and $\lim \limits _{n\rightarrow \infty}\pi(p_n)=\infty$.
  \begin{rk}\label{4..4} Pliss has proved that if $\{p_n(f)\}$ is
  index $i$ stable, then $i\neq 0,d$.\end{rk}

  The following lemma was given by Shaobo Gan, and the proof comes from him
  also.

\begin{lem}\label{4.5}(\cite{1G1})  $f\in C^1(M)$, suppose $\{p_n(f)\}$ is
index $i$ stable, then there exists a subsequence
$\{p_{n_j}\}^\infty_{j=1}$ such that $p_{n_j}$ and $p_{n_{j+1}}$
are homoclinic related.\end{lem}

  Here we just prove the following weaker statement of Gan's
  lemma.

\begin{lem}\label{4.6}( Weaker statement of Gan's lemma)  Suppose
$f\in R$, $\Lambda$ is a non-trivial chain recurrent set of $f$,
$\{p_n(f)\}$ is index $i$ stable and $\lim \limits _{n\rightarrow
\infty}Orb(p_n)=\Lambda$, then there exists a subsequence
$\{p_{n_j}(f)\}^\infty_{j=1}$ such that $p_{n_j}(f)$ and
$p_{n_{j+1}}(f)$ are homoclinic related.\end{lem}

  Before we prove lemma \ref{4.6}, we'll give a few lemmas which will be used in the
  proof.

\begin{lem}\label{4.7} Suppose $A$=$\bigl(\begin{smallmatrix} B&C\\
0&D
\end{smallmatrix}\bigr)$ is a hyperbolic linear map with index
$i\ (i\neq 0,d)$, where $B\in GL(R^i)$ is a contracting map and
$D\in GL(R^{d-i})$ is a expanding map. If there exists $B^\prime
\in GL(R^i)$ an $\varepsilon$-perturbation of $B$ and $B^\prime$
has index different with $i$, then
$A^\prime$=$\bigl(\begin{smallmatrix} B^\prime &C\\ 0&D
\end{smallmatrix}\bigr)$ is an $\varepsilon$-perturbation of $A$
with index different with $i$. In fact, we'll have
$ind(A^\prime)=ind(B^\prime)$.\end{lem}

  With lemma \ref{4.7}, the following lemma is obvious.

\begin{lem}\label{4.8} Suppose $\{\xi^{(n)}\}^\infty_{n=1}$ is
index $i$ stable, then
$\{\xi^{(n)}|_{E^s(\xi^{(n)})}\}^\infty_{n=1}$ is stable
contracting, and at the same time,
$\{\xi^{(n)}|_{E^u(\xi^{(n)})}\}^\infty_{n=1}$ is stable
expanding.\end{lem}

  In \cite{1M3} Ma\~{n}\'{e} has given a necessary condition for bounded stable contracting sequence.

\begin{lem}\label{4.9} (Ma\~{n}\'{e})  If $\{\xi^{(n)}\}^\infty_{n=1}$
is stable contracting and bounded, then there exist
$N_0,l_0,0<\lambda_0<1$ such that if $\pi(\xi^{(n)})>N_0$ we'll
have
 $$   \prod \limits^{[\frac{\pi(\xi_n)}{l_0}]-1}_{j=0}
 \Vert \prod \limits^{l_0-1}_{t=0}\xi^{(n)}_{(jl_0+t)+s}\Vert
 \leq \lambda_0^{[\frac{\pi(\xi^{(n)})}{l_0}]}$$ for any $0\leq
 s<\pi(\xi^{(n)})$.\end{lem}

\noindent{\bf Proof of lemma \ref{4.6}}: Since $\Lambda\subset
P^*_i$ and $f$ is far away from tangency, by proposition 2.1,
$\Lambda$ has an index $i-(l,\lambda)$ dominated splitting
$T|_{\Lambda}M=E\oplus F$. In order to make the proof more
simiplier, here we just suppose $l=1$. Choose a small open
neighborhood $U$ of $\Lambda$, when $U$ is small enough,
$\widetilde{\Lambda}=\bigcap \limits_{j\in
\mathbb{Z}}f^j(\overline{U})$ has an index
$i-(1,\widetilde{\lambda})$ dominated splitting
$T_{\widetilde{\Lambda}}M=\widetilde{E}\oplus \widetilde{F}$ where
$\lambda <\widetilde{\lambda} <1$ and
$\widetilde{E}|_{\Lambda}=E$, $\widetilde{F}|_{\Lambda}=F$.

  Since $\displaystyle{\lim_{n\rightarrow \infty}}Orb(P_n)=\Lambda$,
   we can always suppose $Orb(p_n)\subset \overline{U}$, so $Orb(P_n)\subset \widetilde{\Lambda}$
   and $E^s|_{Orb(p_n)}=\widetilde{E}|_{Orb(p_n)}$, $F^u|_{Orb(p_n)}=\widetilde{F}|_{Orb(p_n)}$.

  By lemma \ref{4.8}, we know that $\{Df|_{E^s(Orb(p_n))}\}^\infty_{n=1}$ is stable contracting and
  $\{Df|_{E^u(Orb(p_n))}\}^\infty_{n=1}$ is stable expanding.
   By lemma \ref{4.9}, there exist$N_0,l_0,0<\lambda_0<1$ such that if $\pi(p_n(f))>N_0$, we have
\begin{equation}\label{1}
  \prod \limits^{[\frac{\pi(p_n)}{l_0}]-1}_{j=0}\Vert Df^{l_0}|_{E^s(f^{jl_0}p_n)}\Vert
  \leq\lambda_0^{[\frac{\pi(p_n)}{l_0}]}
\end{equation}
\begin{equation}\label{2}
 \prod \limits^{[\frac{\pi(p_n)}{l_0}]-1}_{j=0}\Vert
Df^{-l_0}|_{F^u(f^{-jl_0}p_n)}\Vert\leq\lambda_0^{[\frac{\pi(p_n)}{l_0}]}
\end{equation}
  Since $\lim\limits_{n\rightarrow \infty}Orb(p_n)=\Lambda$ and $\Lambda$ is not trivial,
 we have $\lim \limits_{n\rightarrow
\infty}\pi(p_n) \longrightarrow \infty$, then we can always
suppose all the $p_n$ satisfy \eqref{1} and \eqref{2}. For
simplicity, we suppose $l_0=1$ here.

  Choose some $\varepsilon>0$ and $\lambda_1<1$ such that
  $\max\{\widetilde{\lambda},\ \lambda_0\}+\varepsilon<\lambda^2_1<\lambda_1<1$. Now we'll state Pliss lemma in a special context.

\begin{lem}\label{4.10} (Pliss\cite{1Pl}) Given
$0<\lambda_0+\varepsilon<\lambda_1<1$ and $Orb(p_n)\subset
\widetilde{\Lambda}$ such that for some $m\in \mathbb{N}$, we have
$\prod \limits_{j=0}^{t-1}\Vert Df|_{E^s(f^jp_n))}\Vert \leq
(\lambda_0+\varepsilon)^t$ for all $s\geq m$, there exists a
sequence $0\leq n_1<n_2<\cdots$ such that $\prod
\limits_{j=n_r}^{t-1}\Vert Df|_{E^s(f^jp_n))}\Vert \leq
\lambda_1^{t-n_r}$ for all $t\geq n_r$, $r=1,2,\cdots.$\end{lem}

\begin{rk}\label{4.11} The sequence $\{n_j\}^\infty_{j=1}$ we get
above is called the $\lambda_1$-hyperbolic time for bundle
$E^s|_{Orb(p_n)}$. By \eqref{1},\eqref{2}, when $n$ is big enough,
$Orb(p_n)$ will satisfy the assumption of Pliss lemma, so by lemma
\ref{4.10}, there exists $q^+_n\in Orb(p_n)$ such that $\prod
\limits^{t-1}_{j=0}\Vert Df|_{E^s(f^{j}q^+_n)}\Vert\leq
\lambda^t_1$ and $q^-_n\in Orb(p_n)$ such that $\prod
\limits^{t-1}_{j=0}\Vert Df^{-1}|_{F^u(f^{-j}q^-_n)}\Vert\leq
\lambda^t_1$ for all $t>0$.\end{rk}

  Let's denote
  $$S_{n,+}=\{m\in \mathbb{Z}:\prod \limits^{s-1}_{j=0}\Vert Df|_{E^s(f^{m+j}p_n)}\Vert\leq \lambda^s_1\ \text{for\ all}\ s>0\},$$
  $$S_{n,-}=\{m\in \mathbb{Z}:\prod \limits^{s-1}_{j=0}\Vert
Df^{-1}|_{F^u(f^{m-j}p_n)}\Vert\leq \lambda^s_1\ \text{for\ all}\
s>0\}.$$ Then $S_{n,+}$ is the set of $\lambda_1$ hyperbolic time
for bundle $E^s|_{Orb(p_n)}$ and $S_{n,-}$ is the set of
hyperbolic time for bundle $F^u|_{Orb(p_n)}$. From remark 4.11,
the set $S_{n,+}$ and $S_{n,-}$ are not empty. We denote
$S_n=S_{n,+}\bigcap S_{n,-}$.

\begin{lem}\label{4.12} $S_n\neq \phi$.\end{lem}

 \noindent{\bf Proof}$\colon$ Here for $a,b\in \mathbb{Z}$ and $a<b$, we denote
$(a,b)_\mathbb{Z}=\{c|\ c\in \mathbb{Z}\  and\  a<c<b\}$.

  Now suppose the lemma is false, we can choose $\{b_{n,s},\ b_{n,s+1}\}\subset S_{n,-}$ such that we have
   $b_{n,s+1}>b_{n,s}$, $(b_{n,s}, b_{n,s+1})_{\mathbb{Z}}\bigcap S_{n,-}= \phi$ and
   $a_{n,t}\in (b_{n,s},b_{n,s+1})_{\mathbb{Z}}\bigcap S_{n,+}$, then $b_{n,s},b_{n,s+1}\notin S_{n,+}$.

  We claim that for $0<k\leq b_{n,s+1}-b_{n,s}-1$, we have
  $\prod \limits^{k-1}_{j=0}\Vert Df^{-1}|_{F^u(f^{b_{n,s}+j+1}p_n)}\Vert \geq \lambda_1^k.$

  \noindent{\bf Proof of the claim}:  We'll use induction to give a proof.

When $k=1$, since $b_{n,s}+1\notin S_{n,-}$, we have $\Vert
Df^{-1}|_{F^u(f^{b_{n,s}+1}p_n)}\Vert>\lambda_1$.

Now suppose the claim is true for all $1\leq k \leq k_0-1$ where $1<k_0\leq b_{n,s+1}-b_{n,s}-1$,
  and we suppose the claim is false for $k_0$, it means that
\begin{equation}\label{3}
\prod \limits^{k_0-1}_{j=0}\Vert
Df^{-1}|_{F^u(f^{b_{n,s}+j+1}p_n)}\Vert \leq \lambda_1^{k_0}.
\end{equation} Then by the assumption above that the claim is true for $1\leq k\leq k_0-1$, we
have \begin{equation}\label{4}
  \prod \limits^{k-1}_{j=0}\Vert
  Df^{-1}|_{F^u(f^{b_{n,s+j+1}}p_n)}\Vert \geq \lambda_1^k
     \end{equation}
  By \eqref{3} and \eqref{4}, we get that
  $\prod \limits^{k_0-1}_{j=k}\Vert Df^{-1}|_{F^{u}(f^{b_{n,s}+j+1}p_n)}\Vert < \lambda_1^{k_0-k}$
  for $1\leq k\leq k_0-1$. It's equivalent to say that
  \begin{equation}\label{5}
\prod \limits^{k-1}_{j=0}\Vert
Df^{-1}|_{F^{u}(f^{b_{n,s}+k_0-j}p_n)}\Vert < \lambda_1^{k}
\quad \text{for} \ 1\leq k\leq k_0-1
  \end{equation}
  By \eqref{3} and \eqref{5}, we get that
  \begin{equation}\label{6}
\prod \limits^{k-1}_{j=0}\Vert
Df^{-1}|_{F^{u}(f^{b_{n,s}+k_0-j}p_n)}\Vert \leq \lambda_1^{k}
\quad \text{for} \ 1\leq k\leq k_0
  \end{equation}
  When $k>k_0$,by \eqref{6} and the fact $b_{n,s}\in S_{n,-}$, we have
$$
  \prod \limits^{k-1}_{j=0}\Vert Df^{-1}|_{F^{u}(f^{b_{n,s}+k_0-j})}\Vert
   =\prod \limits^{k_0-1}_{j=0}\Vert Df^{-1}|_{F^{u}(f^{b_{n,s}+k_0-j})}\Vert\cdot
  \prod \limits^{k-k_0-1}_{j=0}\Vert Df^{-1}|_{F^{u}(f^{b_{n,s}-j})}\Vert
   <\lambda^{k_0}_1\cdot \lambda_1^{k-k_0}=\lambda_1^{k},$$
    it means $b_{n,s}+k_0\in S_{n,-}$, it's a contradiction since $b_{n,s}+k_0\in(b_{n,s},b_{n,s+1})_{\mathbb{Z}}$, so we finish
the induction.                           \qed

  By the claim above, for $0<k\leq b_{n,s+1}-b_{n,s}-1$, we have
  \begin{equation}\label{7}
\prod \limits^{k-1}_{j=0}\Vert
Df^{-1}|_{F^u(f^{b_{n,s}+j+1}p_n)}\Vert \geq\lambda_1^k.
  \end{equation}
  Since on $\widetilde{\Lambda}$, $\widetilde{E}\oplus \widetilde{F}$
   is an index $i-(1,\widetilde{\lambda})$ dominated splitting, we have
  $$\prod \limits_{j=0}^{k-1}(\Vert Df|_{\widetilde{E}(f^{b_{n,s}+j}p_n)}\Vert \cdot
  \Vert Df^{-1}|_{\widetilde{F}(f^{b_{n,s}+j+1}p_n)}\Vert) <\widetilde{\lambda}^k.$$
  By \eqref{7} and
  $\widetilde{E}|_{Orb(p_n)}=E^s|_{Orb(p_n)}$,$\widetilde{F}|_{Orb(p_n)}=F^u|_{Orb(p_n)}$,
we'll get
\begin{equation}\label{8}
\prod \limits_{j=0}^{k-1}\Vert Df|_{E^s(f^{b_{n,s}+j}p_n)}\Vert <
\frac{\widetilde{\lambda}^k}{\lambda^k_1} \underset{(\tiny
\widetilde{\lambda}<\lambda_1^2<1)}{<}\lambda_1^k \quad \text{for}\
1<k\leq b_{n,s+1}-b_{n,s}-1.
\end{equation}
  When $k>b_{n,s+1}-b_{n,s}-1$, let $k=(a_{n,t}-b_{n,s})+(k-a_{n,t})$, by \eqref{8} and
  $a_{n,t}\in S_{n,+}$,
  \begin{eqnarray}\label{9}
\nonumber    \prod \limits_{j=0}^{k-1}\Vert
Df|_{E^s(f^{b_{n,s}+j}p_n)}\Vert &=& \prod
\limits_{j=0}^{a_{n,t}-b_{n,s}-1}\Vert
Df|_{E^s(f^{b_{n,s}+j}p_n)}\Vert\cdot \prod
\limits_{j=0}^{k-a_{n,t}-1}\Vert Df|_{E^s(f^{a_{n,t}+j}p_n)}\Vert \\
     &<& \lambda^{a_{n,t}-b_{n,s}}_1\cdot
\lambda_1^{k-a_{n,t}}=\lambda_1^{k-b_{n,s}}
  \end{eqnarray}
  By \eqref{8} and \eqref{9}, we get $b_{n,s}\in S_{n,+}$, so $S_{n,+}\bigcap S_{n,-}\neq \phi$,
  it's a contradiction with our assumption, so we finish the proof of lemma 4.12. \qed

  Now let's continue the proof of lemma \ref{4.6}, we need the following two lemmas to show that for $a_{n}\in S_n$,
  the point $f^{a_n}(p_n)$ will have uniform size of stable manifold and unstable manifold.

  Let $I_1=(-1,1)^i$ and $I_\varepsilon=(-\varepsilon,\ \varepsilon)^i$,
  denote by $Emb^1(I,M)$ the set of $C^1$-embedding of $I_1$ on $M$, recall by
  \cite{1HPS}
  that $\widetilde{\Lambda}$ has a dominated splitting $\widetilde{E}\oplus\widetilde{F}$ implies the following.

\begin{lem}\label{4.13} There exist two continuous function
$\Phi^{cs}:\ \widetilde{\Lambda}\longrightarrow Emb^1(I,M)$ and
$\Phi^{cu}:\ \widetilde{\Lambda}\longrightarrow Emb^1(I,M)$ such
that, with $W^{cs}_\varepsilon (x)=\Phi^{cs}(x)I_\varepsilon$ and
$W^{cu}_\varepsilon (x)=\Phi^{cu}(x)I_\varepsilon$, the following
properties hold:
\begin{itemize}
  \item[a)] $T_xW^{cs}_\varepsilon=\widetilde{E}(x)$ and $T_xW^{cu}_\varepsilon=\widetilde{F}(x)$,
  \item[b)]  For all $0<\varepsilon_1<1$, there exists $\varepsilon_2$ such that
  $f(W^{cs}_{\varepsilon_2}(x))\subset W^{cs}_{\varepsilon_1}(f(x))$ and
  $f^{-1}(W^{cu}_{\varepsilon_2}(x))\subset W^{cu}_{\varepsilon_1}(f^{-1}(x))$.
  \item[c)] For all $0<\varepsilon<1$, there exists $\delta>0$ such that if $y_1,\ y_2\in\widetilde{\Lambda}$ and
  $d(y_1,\ y_2)<\delta$, then $W^{cs}_{\varepsilon}(y_1)\h W^{cu}_{\varepsilon}(y_2)\neq \phi$.
\end{itemize}\end{lem}

\begin{cor} \label{4.14} (\cite{1PS1}) For any $0<\lambda<1$, there
exists $\varepsilon>0$ such that for $x\in \widetilde{\Lambda}$
which satisfies $\prod \limits_{j=0}^{n-1}\Vert
Df|_{\widetilde{E}(f^jx)}\Vert \leq \lambda^n$ for all $n>0$, then
$diam(f^n(W^{cs}_{\varepsilon}))\longrightarrow 0$, i.e. the
central stable manifold of $x$ with size $\varepsilon$ is in fact
a stable manifold.\end{cor}

  Now for $\lambda_1$, using corollary 4.14, we can get an $\varepsilon>0$. It means that for any $a_n\in S_n$,
  denote $q_n=f^{a_n}(p_n)$, then $W^{cs}_{\varepsilon}(q_n)$ is a stable manifold and
  $W^{cu}_{\varepsilon}(q_n)$ is an unstable manifold.
  For this $\varepsilon>0$, use c) of lemma \ref{4.13}, we can fix a $\delta$.
  Choose a subsequence $\{n_i\}$ such that $d(q_{n_i},q_{n_{i+1}})\leq \delta$,
  then by c) of lemma \ref{4.13}, we know $W^{cu}_{\varepsilon}(q_{n_i})\h W^{cs}_{\varepsilon}(q_{n_{i+1}})\neq \phi$
  and $W^{cu}_{\varepsilon}(q_{n_{i+1}})\h W^{cs}_{\varepsilon}(q_{n_{i}})\neq \phi$.
  Since the local central stable manifold and local central unstable manifold of $q_{n_i}$ have
dynamical meaning, we know that $Orb(q_{n_i})$ and
$Orb(q_{n_{i+1}})$ are homoclinic related.
  \qed

\begin{rk} \label{4.15}  In the proof of lemma \ref{4.6} we suppose the
set $\Lambda$ has 1-step dominated splitting, that means $l=1$,
and we suppose $l_0=1$ there also, they are just in order to make
the proof more simplier. In the rest part of the paper, usually we
don't use such assumption any more, if we use it we'll point
out.\end{rk}

  Now let's consider a sequence of periodic points which are not index stable.

\begin{lem} \label{4.16} Suppose $f\in R$, $\lim
\limits_{n\rightarrow \infty}g_n=f$, $\{p_n(g_n)\}^{\infty}_{n=1}$
is a family of index $i$ periodic points ($i\neq \ 0,\ d$) and
$\lim \limits_{n\rightarrow \infty}\pi(p_n)\longrightarrow
\infty$. If there exist $\lambda_n \longrightarrow 1^{-}$ and
$\lim \limits_{n\rightarrow \infty}l_n\longrightarrow \infty$ such
that $\lim \limits_{n\rightarrow
\infty}\frac{\pi(p_n)}{l_n}\longrightarrow \infty$ and $\prod
\limits^{[\frac{\pi(p_n)}{l_n}]-1}_{j=0}\Vert
Dg_n^{l_n}|_{E^s(g_n^{jl_n}(p_n))}\Vert \geq
\lambda_n^{[\frac{\pi(p_n)}{l_n}]}$, then for any $\varepsilon>0$
and $N>0$, there exists an $n_0>N$ and $g_{n_0}^\prime$ is an
$\varepsilon$-perturbation of $g_{n_0}$ such that
$p_{n_0}(g_{n_0})$ is an index $i-1$ periodic point of
$g_n^\prime$.\end{lem}

\noindent{\bf Proof}: Fix $N$, consider the periodic sequence of
linear maps $\{\xi^n:\ \xi^n=Dg_n|_{E^s(Orb(p_n))}\}_{n\geq N}$, they are
all contracting maps. We claim that $\{ \xi^n \}$ are not stable
contracting.\\

  \noindent{\bf Proof of the claim}: If $\{ \xi^n \}$ is stable contracting,
  by lemma \ref{4.9}, there exist $N_0,\ l_0,\ 0<\lambda_0<1$ such that if $\pi(\xi^n)>N_0$,
  we have \begin{equation}\label{10}
\prod \limits^{[\frac{\pi(p_n)}{l_0}]-1}_{j=0}\Vert
Dg_n^{l_0}|_{E^s(g_n^{jl_0}p_n)}\Vert\leq
\lambda_0^{[\frac{\pi(p_n)}{l_0}]}
  \end{equation}
Choose some $N_1$ big enough such that for $n\geq N_1$, we have
$\lambda_n\geq \ \lambda^*>\ \lambda_0$ for some $\lambda^*\in \
(\lambda_0,1)$, then by $\lim \limits_{n\rightarrow
\infty}\frac{\pi(p_n)}{l_n}\longrightarrow \infty$ and $\lim
\limits_{n\rightarrow \infty}l_n \longrightarrow \infty$, when $n$
is big enough, we have $\pi(p_n)\gg \ l_n\gg \ \max\{l_0,N_0\}$
and from $\prod \limits^{[\frac{\pi(p_n)}{l_0}]-1}_{j=0}\Vert
Dg_n^{l_n}| _{E^s(g_ n^{jl_n}p_n)}\Vert\geq
\lambda_n^{[\frac{\pi(p_n)}{l_n}]}>(\lambda^*)^{[\frac{\pi(p_n)}{l_n}]}$,
we'll get $\prod \limits^{[\frac{\pi(p_n)}{l_0}]-1}_{j=0}\Vert
Dg_n^{l_0}| _{E^s(g_ n^{jl_0}p_n)}\Vert \geq
\lambda_0^{[\frac{\pi(p_n)}{l_n}]}>\lambda_0^{[\frac{\pi(p_n)}{l_0}]}$,
It's a contradiction with \eqref{10}. \qed

Since $\{\xi^n\}_{n\geq N}$ isn't stable contracting, for
$\varepsilon>0$, there exists a sequence $\{n_i\}$ and
$\{\eta^{n_i}\}$ such that $\eta^{n_i}$ is an
$\varepsilon$-perturbation of $\xi^{n_i}$ and $\eta^{n_i}$ has
index smaller than $i$. Since $\{\xi^{n_i}\}$ is bounded and $\lim
\limits_{n\rightarrow \infty}\pi(p_n)\longrightarrow \infty$, by
\cite{1BGV}'s work, for $n_i$ big enough, we can in fact get
$\eta^{n_i}$ with index $i-1$. By lemma \ref{4.7}, there exists
$\{A|_{Orb(p_n)}\}_{n\geq 0}$ an $\varepsilon$-perturbation of
$\{Dg_n|_{Orb(p_n)}\}$ such that $\{A|_{Orb(p_n)}\}$ has index
$i-1$.
  Now we need the following version of Franks lemma.

\begin{lem} \label{Franks}(Franks lemma) Suppose $p_n$ is a periodic point of
$g_n$, $A|_{Orb(p_n)}$ is an $\varepsilon$-perturbation of
$\{Dg_n|_{Orb(p_n)}\}$, then for any neighborhood $U$ of
$Orb(p_n)$, there exists $g_n^\prime$ such that $g_n^\prime\equiv
g_n$ on $(M\setminus U)\bigcup Orb(p_n)$,
$d_{C^1}(g_n,g_n^\prime)<\varepsilon$ and
$\{Dg_n^\prime|_{orb(p_n)}\}=\{A|_{Orb(p_n)}\}$.\end{lem}

From Franks lemma, we can change the derivative map along
$T_{Orb(p_{n_i})}M$ to be $\{A|_{Orb(p_n)}\}$ and get a new map
$g_{n_i}^\prime$ such that $p_{n_i}(g_{n_i})$ is index $i-1$
periodic point of $g_{n_i}^\prime$. \qed

\subsection{Weakly selecting lemma}
  Liao's selecting lemma is a powerful shadowing lemma for
  non-uniformly hyperbolic system, with it, we can not only get
 a lot of periodic points like what the standard shadowing lemma
 can do, we can even let the periodic points have hyperbolic property
 as weak as we like. Liao at first used this lemma to study minimal
 non-hyperbolic set and proved the $\Omega$-stable conjecture for diffeomorphisms
  in dimension 2 and for flow without singularity in dimension 3.
  \cite{1G2} \cite{1G3} \cite{1GW2} \cite{1W0} use the same idea proved structure($\Omega$) stability conjecture for flows
  without singularity in any dimension. Until now, the most important
  papers about selecting lemma are \cite{1GW1},\cite{1W3}, \cite{1W4} and there contain more
  details about selecting lemma.

In this subsection and the next, we'll show what will happen if
all the conditions in weakly selecting lemma are satisfied. The
main result in this subsection is lemma \ref{4.21} (The weakly
selecting lemma). Now let's state the selecting lemma at first.

\begin{pro}\label{4.18}(Liao) Let $\Lambda$ be a compact invariant
set of $f$ with index $i-(l,\lambda)$ dominated splitting
$E^{cs}\oplus F^{cu}$. Assume that
\begin{itemize}
  \item[a)] there is a point $b\in \Lambda$ satisfying $\prod \limits_{j=0}^{n-1}\Vert Df^l|_{E^{cs}(f^{jl}b)}\Vert\geq 1$
for all $n\geq 1$.
  \item[b)] (The tilda condition) there are $\lambda_1$ and $\lambda_2$ with
  $\lambda<\lambda_1<\lambda_2<1$
such that for any $x\in \Lambda$ satisfying $\prod
\limits_{j=0}^{n-1}\Vert Df^l|_{E^{cs}(f^{jl}x)}\Vert\geq
{\lambda_2}^n$ for all $n\geq 1$, $\omega(x)$ contains a point
$c\in \Lambda$ satisfying $\prod \limits_{j=0}^{n-1}\Vert
Df^l|_{E^{cs}(f^{jl}c)}\Vert\leq \lambda_1^n$ for all $n\geq
1$.\end{itemize} Then for any $\lambda_3$ and $\lambda_4$ with
$\lambda_2<\lambda_3<\lambda_4<1$ and any neighborhood $U$ of
$\Lambda$, there exists a hyperbolic periodic orbit $Orb(q)$ of
$f$ of index $i$ contained entirely in $U$ with a point $q\in
Orb(q)$ such that
\begin{equation}\label{*}
\prod \limits^{m-1}_{j=0}\Vert Df^l|_{E^{cs}(f^{jl}q)} \Vert \leq
\lambda_4^m,\quad for\ m=1,\cdots,\ \pi_l(q)
\end{equation}
\begin{equation}\label{**}
\prod \limits^{\pi_l(q)-1}_{j=\pi_l(q)-m}\Vert
Df^l|_{E^{cs}(f^{jl}q)}\Vert \geq \lambda_3^{m} \quad for \
m=1,\cdots,\ \pi_l(q)
\end{equation}
where $\pi_l(q)$ is the period of $q$ for the map $f^l$. The
similar assertion for $F^{cu}$ holds respecting $f^{-1}$.\end{pro}

\begin{rk} \label{4.19} It's easy to know $\pi(q)\geq \pi_l(q)$. Since
$f^{l\cdot \pi_l(q)}(q)=q$, it's obvious that \eqref{*} and
\eqref{**} are true for all $m\in \mathbb{N}$. In the selecting
lemma, when $\lambda_3$ and $\lambda_4$ are fixed, we can indeed
find a sequence of periodic points $\{q_n\}$ satisfying \eqref{*}
and \eqref{**} and $\overline{\lim \limits_{n\rightarrow
\infty}}Orb(q_n)\subset \Lambda$. If $f$ is a Kupuka-Smale
diffeomorphism, especially when $f\in R$, we can let $\lim
\limits_{n\rightarrow \infty}\pi_l(q_n)\longrightarrow \infty$,
then we'll have $\lim \limits_{n\rightarrow
\infty}\pi(q_n)\longrightarrow \infty$ at the same time.\end{rk}

\begin{cor} \label{4.20} $f\in R$, let $\Lambda$ be a compact chain
recurrent set of $f$ with index $i-(l_0,\lambda)$ dominated
splitting $E^{cs}\oplus F^{cu}$ ($1\leq\ i\leq\ d-1$). Assume that
the splitting satisfies all the conditions of selecting lemma for
all $l_n=nl_0$ ($n\in \mathbb{N}$) but with the same parameters
$\lambda\ <\lambda_1\ <\lambda_2\ <1$, then for any sequence
$\{(\lambda_{n,3},\ \lambda_{n,4})\}^\infty_{n=1}$ satisfying
$\lambda_2<\lambda_{1,3}<\lambda_{1,4}<\lambda_{2,3}<\lambda_{2,3}<\cdots$
where $\lambda_{n,3}\longrightarrow 1^{-}$, there exists a family
of periodic points $\{q_n(f)\}$ with index $i$ such that
\begin{itemize}
\item[a)] $\lim \limits_{n\rightarrow \infty}
\pi_{l_n}(q_n(f))\longrightarrow \infty$.

\item[b)]
\begin{subequations} \begin{align}
\prod \limits^{m-1}_{j=0}\Vert Df^{l_n}|_{E^s(f^{jl_n}q_n)}\Vert &\leq \lambda^m_{n,4}\label{11}\\
\prod \limits^{\pi_{l_n}(q_n)-1}_{j=\pi_{l_n}(q_n)-m}\Vert
Df^{l_n}|_{E^s(f^{jl_n}q_n)}\Vert &\geq \lambda^m_{n,3}\quad
\text{for}\ m\in \mathbb{N} \label{12}
\end{align}
\end{subequations}

\item[c)] $\overline{\lim \limits_{n\rightarrow
\infty}}Orb(q_n)\subset \Lambda$.

\item[d)] $\Lambda \subset H(q_n(f))$\quad for all $n$.
\end{itemize}
\end{cor}

\Pf: At first, let's fix $\lambda_2<\lambda_{1,3}<\lambda_{1,4}<1$
and a small neighborhood $U$ of $\Lambda$ small enough such that
the maximal invariant set $\widetilde{\Lambda}$ of $\overline{U}$
has index $i-(l_0,\widetilde{\lambda})$ dominated splitting with
$\widetilde{\lambda}<\lambda_2$, we denote the dominated splitting
still by $E^{cs}_i\oplus F^{cu}_{i+1}$. (If $q$ is an index $i$
periodic point in $\widetilde{\Lambda}$, then we denote
$E^{cs}_i\oplus F^{cu}_{i+1}|_{Orb(q)}=E^s\oplus F^u|_{Orb(q)}$).
Now using selecting lemma, with remark \ref{4.19}, we can find a
family of periodic points $\{q_{1,m}(f)\}^\infty_{m=1}$ with index
$i$ satisfying b), $\overline{\lim \limits_{n\rightarrow
\infty}}(q_{1,m})\subset \Lambda$, $\lim \limits_{m\rightarrow
\infty} \pi_{l_1}(q_{1,n})\longrightarrow \infty$ and
$Orb(q_{1,m}(f))\subset \widetilde{\Lambda}$.

  Since $\widetilde{\Lambda}$ has an index $i-(l_1,\widetilde{\lambda})$ dominated splitting
$E^{cs}_{\widetilde{\Lambda}}\oplus F^{cu}_{\widetilde{\Lambda}}$,
from \eqref{12} we can know \begin{equation*} \prod
\limits_{j=\pi_{l_1}(q_{1,m})-t+1}^{\pi_{l_1}(q_{1,m})}\Vert
Df^{-l_1}|_{F^{cu}(f^{jl_1}q_{1,m})}\Vert \leq
\widetilde{\lambda}^t/ \prod
\limits_{j=\pi_{l_1}(q_{1,m})-t}^{\pi_{l_1}(q_{1,m})-1}\Vert
Df^{l_1}|_{E^{cs}(f^{jl_1}q_{1,m})}\Vert\leq
(\frac{\widetilde{\lambda}}{\lambda_{1,3}})^t\ \text{for}\ _{(t\in
\mathbb{N})},
\end{equation*}
it equivalent with \begin{equation}\label{13} \prod
\limits_{t=0}^{m-1}\Vert
Df^{-l_1}|_{F^{cu}(f^{-jl_1}q_{1,m})}\Vert \leq
(\frac{\widetilde{\lambda}}{\lambda_{1,3}})^t \quad
  \text{for}\ t\in \mathbb{N}.\end{equation}

  From \eqref{11}, \eqref{12}, by lemma \ref{4.13}, Corollary \ref{4.14} and
  $\frac{\widetilde{\lambda}}{\lambda_{1,3}}<1$, we can know that for some $\varepsilon_1$,
  $q_{1,n}$ will have uniformly size of stable manifold $W_{\varepsilon_1}^{s}(q_{1,n})$
  and uniform size of unstable manifold $W_{\varepsilon_1}^{u}(q_{1,n})$
  and there exists a subsequence $\{q_{1,n_{j}}\}_{j=1}^{\infty}$ such
   that they are homoclinic related with each other, so $H(q_{1,n_1})=H(q_{1,n_2})=\cdots$, with
   $\overline{\lim \limits_{j\rightarrow \infty}}Orb(q_{1,n_j})\subset \Lambda$,
   we know $\Lambda\bigcap H(q_{1,n_j})\neq \phi$. Since $f\in R$, $H(q_{1,n_j})$ should be a chain recurrent class.
    Because $\Lambda$ is a chain recurrent set, we have $\Lambda\subset H(q_{1,n_j})$,
    let $q_1=q_{1,n_j}$ for some $j$ big enough, then $q_1$ satisfies $a)$, $d)$.

Now consider $0<\lambda_2<\lambda_{2,3}<\lambda_{2,4}<1$,
$E^{cs}_{\Lambda}\oplus F^{cu}_{\Lambda}$ is obviously an index
$i-(l_2,\lambda)$ dominated splitting of $\Lambda$ and by the
assumption, the splitting satisfy the conditions of selecting
lemma for $l_2,\ \lambda<\lambda_1<\lambda_2<1$, so repeat the
above argument, we can get a family of periodic points
$\{q_{2,n}(f)\}_{n=1}^{\infty}$ satisfying $b)$, $d)$,
$\overline{\lim \limits_{n\rightarrow \infty}}Orb(q_{2,n})\subset
\Lambda$, $\Lambda\subset H(q_{2,1},f)=\cdots=H(q_{2,n},f)=\cdots$
and $\lim \limits_{n\rightarrow
\infty}\pi_{l_2}(q_{2,n}(f))\longrightarrow \infty$. When $n_0$ is
big enough, we'll have $\pi_{l_2}(q_{2,n_0})>\pi_{l_1}(q_1)$ and
$Orb(q_{2,n_0})$ is near $\Lambda$ more than $Orb(q_1)$. Let
$q_2=q_{2,n_0}$, continue the above argument for $l_n$ and
$\lambda_2<\lambda_{n,3}<\lambda_{n,4}<1$, we can get
$\{q_n\}^{\infty}_{n=1}$ which we need.
      \qed

  The following weakly selecting lemma shows when the conditions of the above corollary will be satisfied.

\begin{lem}\label{4.21}(Weakly selecting lemma) Let $f\in R$, $\Lambda$ be
a compact invariant set of $f$ with index $i-(l_0,\ \lambda)$
dominated splitting $E^{cs}\oplus F^{cu}$ $(1\leq i\leq d-1)$.
Assume that
\begin{itemize}

  \item[a)](Non-hyperbolic condition) the bundle $E^{cs}$ is not contracting,

  \item[b)](Strong tilda condition) there are $\lambda_2<1$ and $l_0^\prime>1$
  such that for any $x\in \Lambda$, $\omega(x)$ contains a point $c\in \Lambda$ satisfying
  $\prod \limits_{j=0}^{n-1}\Vert Df^{l_0^\prime}|_{E^{cs}(f^{jl_0^\prime}c)}\Vert\leq \lambda_2^n$ for all $n\geq 1$.
\end{itemize}

  Then for any $l_n =n\cdot(l_0\cdot l_0^\prime)$ and any sequence
  $\{(\lambda_{n,3},\lambda_{n,4})\}^\infty_{n=1}$ satisfying $max \{\lambda^{l_0^\prime}, \lambda_2\}
  <\lambda_{1,3}<\lambda_{1,4}<\cdots<\lambda_{n,3}<\lambda_{n,4}<\cdots$ where
  $\lambda_{n,3}\longrightarrow 1^-$, there exists a family of periodic points $\{q_n(f)\}$ with index $i$ such that
\begin{itemize}

  \item[$\bullet$] $\lim \limits_{n\rightarrow \infty} \pi_{l_n}(q_n(f))\longrightarrow \infty$

  \item[$\bullet$] $\prod \limits_{j=0}^{m-1}\Vert Df^{l_n}|_{E^s(f^{jl_n}q_n)}\Vert \leq
  \lambda_{n,4}^m$ and $\prod \limits_{j=\pi_{l_n}(q_n)-m}^{\pi_{l_n}(q_n)-1}\Vert Df^{l_n}|_{E^s(f^{jl_n}q_n)}\Vert\geq
   \lambda_{n,3}^m$ for $m\geq 1$

  \item[$\bullet$] $\overline{\lim \limits_{n\rightarrow \infty}}Orb(q_n)\subset \Lambda$

  \item[$\bullet$] $\Lambda \subset H(q_n(f))$ for $n \geq 1$.
\end{itemize}
\end{lem}

\Pf Since $E^{cs}_\Lambda\oplus F^{cu}_\Lambda$ is a
$(l_0,\lambda)$ dominated splitting and $l_1=l_0\cdot l_0^\prime$,
it should be a $(l_1,\lambda^{l_0^\prime})$ dominated splitting
also. Choose $\lambda_2^\prime,\ \lambda_1$ such that $\max
\{\lambda^{l_0^\prime},\
\lambda_2\}<\lambda_1<\lambda_2^\prime<\lambda_{1,3}$, we'll show
that the splitting $E^{cs}_\Lambda\oplus F^{cu}_\Lambda$ and the
$l_1,\ \lambda^{l_0^\prime}<\lambda_1<\lambda_2^\prime<1$ will
satisfy all conditions of corollary 4.20, equivalent, we'll show
the splitting $E^{cs}_\Lambda\oplus F^{cu}_\Lambda$, $l_n$ and
$\lambda^{l_0^\prime}<\lambda_1<\lambda_2^\prime<1$ will satisfy
the condition of selecting lemma for all $n\geq 1$.
\begin{itemize}
\item[0)] Since $E^{cs}_\Lambda\oplus F^{cu}_\Lambda$ is a $(l_1,\
\lambda^{l_0^\prime})$ dominated splitting and $l_n=n\cdot l_1$,
$E^{cs}_\Lambda\oplus F^{cu}_\Lambda$ is a $(l_n,\
\lambda^{l_0^\prime})$ dominated splitting also.

\item[1)]  Here we need the following lemma:

\begin{lem}\label{4.22} Let $\Lambda$ be a compact invariant set of $f$,
$E^{cs}_\Lambda$ is an continuous invariant bundle on $\Lambda$,
and $dim(E^{cs}(x))=i$ for any $x\in \Lambda$ where $i\neq 0$,
suppose $l\in \mathbb{N}$, if for any $x\in \Lambda$, there exists
an $n$ such that $\prod \limits_{j=0}^{n-1}\Vert
Df^l|_{E^{cs}(f^{jl}x)}\Vert< 1$, then $E^{cs}_\Lambda$ is a
contracting bundle.\end{lem}

Since we know $E^{cs}_\Lambda$ is continuous but not contracting,
so for any $l_n$, there exists $b_n$, such that $\prod
\limits_{j=0}^{n-1}\Vert Df^l_n|_{E^{cs}(f^{jl_n}b_n)}\Vert\geq 1$
for all $m\geq 1$.

\item[2)]  For any $x\in \Lambda$, $\omega(x)$ contains a point
$c_n\in \Lambda$ such that $\prod \limits_{j=0}^{nl_0m-1}\Vert
Df^{l_0^\prime}|_{E^{cs}(f^{jl_0^\prime}c_n)}\Vert\leq
\lambda_2^{nl_0m}$ for all $m\geq 1$, since
$$\prod \limits_{j=0}^{nl_0m-1}\Vert Df^{l_0^\prime}|_{E^{cs}(f^{jl_0^\prime}c_n)}\Vert\geq
\prod \limits_{j=0}^{m-1}\Vert
Df^{nl_0l_0^\prime}|_{E^{cs}(f^{jnl_0l_0^\prime}c_n)}\Vert=\prod
\limits_{j=0}^{m-1}\Vert Df^{l_n}|_{E^{cs}(f^{jl_n}c_n)}\Vert ,$$
we have that $\prod \limits_{j=0}^{m-1}\Vert
Df^{l_n}|_{E^{cs}(f^{jl_n}c_n)}\Vert \leq \lambda_2^{mnl_0}\leq
\lambda_2^m$ for all $m\geq 1$.
\end{itemize}

\begin{rk}\label{4.23}  In $b)$ of weakly selecting lemma, we don't give any restriction on $x$, so
$b)$ is in fact more stronger than the tilda condition, that's why
we call the condition $b)$ in weakly selecting lemma the strong
tilda condition.\end{rk}
  By 0), 1), 2) above and corollary \ref{4.20}, we proved the lemma.                \qed
\subsection{Transition}

Transition was introduced in \cite{1BDP} at first, there they
consider a special linear system with a special property called
transition and use it to study homoclinic class. Here I prefer to
use a little different way to state it, the notation and
definition are basically copy from \cite{1BDP}. The main result in
this subsection is corollary \ref{4.26}. We begin by giving some
definitions.

Given a set $\cal A$, a word with letters in $\cal A$ is a finite
sequence of $\cal A$, its length is the number of letters
composing it. The set of words admits a natural semi-group
structure: the product of the word $[a]=(a_1,\cdots,a_n)$ by
$[b]=(b_1,\cdots,b_l)$ is
$[a]\cdot[b]=(a_1,\cdots,a_n,b_1,\cdots,b_l)$. We say that a word
$[a]$ is not a power if $[a]\neq[b]^k$ for every word $[b]$ and
$k>1$.

Here we'll use some special words. Let's fix $f\in C^1(M)$, for
any $x\in Per(f)$, we write $[x]=(f^{\pi(x)-1}(x)),\cdots,x)$ and
$\{x\}=(Df(f^{\pi(x)-1}(x)),\cdots,Df(x))$. We call a word
$[a]=(a_k,\cdots,a_1)$ with letters in $M$ is a finite
$\varepsilon$-pseudo orbit if $d(f(a_i),a_{i+1})\leq \varepsilon$
for $1\leq i\leq k-1$, if $\varepsilon=0$, that means
$f(a_i)=a_{i+1}$ for $1\leq i\leq k-1$, then we call $[a]$ is a
finite segment of orbit. We always denote
$\{a\}=(Df(a_k),\cdots,Df(a_1))$.

Suppose we have a finite orbit $[a]=(a_n,\cdots,a_1)$ and an
$\varepsilon$-pseudo orbit $[b]=(b_l,\cdots,b_1)$, we say $[b]$ is
$\delta$-shadowed by $[a]$ if $n=l$ and $d(a_i,b_i)\leq
\varepsilon$ for $1\leq i\leq n$. We say $\{a\}$ is $\delta$-close
to $\{b\}$ if $n=l$ and $\Vert Df(a_i)-Df(b_i)\Vert \leq \delta$
for $1\leq i \leq n$.

Suppose $H(p,f)$ is a non-trivial homoclinic class, we say
$H(p,f)$ has $\varepsilon$-transition property if : for any finite
hyperbolic periodic points ${p_1,\cdots,p_n}$ in $H(p,f)$ which
are homoclinic related with each other, there exist finite orbits
$[t^{i,j}]=(t^{i,j}_{k(i,j)},\cdots,t^{i,j}_1)$ for any $(i,j)\in
\{1,\cdots,n\}^2$ where $k(i,j)$ is the length of $[t^{i,j}]$,
such that, for every $m\in \mathbb{N}$, $l=(i_1,\cdots,i_m) \in
\{1,\cdots,n\}^m$ and $\alpha=(\alpha_1,\cdots,\alpha_m)\in
\mathbb{N}^m$ where the word
$(({i_1},\alpha_1),\cdots,({i_m},\alpha_m))$ with letters in
$\mathbb{N}\times \mathbb{N}$ is not a power, the pseudo orbit
$[w(l,\alpha)]=[t^{i_m,i_1}]\cdot [p_{i_m}]^{\alpha_m}\cdot
[t^{i_{m-1},i_m}]\cdot [p_{i_{m-1}}]^{\alpha_{m-1}}\cdot\cdots
\cdot[t^{i_1,i_2}]\cdot [p_{i_1}]^{\alpha_1}$ is an
$\varepsilon$-pseudo orbit and there is a periodic orbit
$Orb(q(l,\alpha))\subset H(p,f)$ such that:
\begin{itemize}

\item[a)] the length of $[w(l,\alpha)]$ is $\pi(q(l,\alpha))$ and
$[q(l,\alpha)]$ $\varepsilon$-shadow the pseudo orbit
$[w(l,\alpha)]$.

\item[b)] the word $\{q(l,\alpha)\}$ is $\varepsilon$-close to
$\{w(l,\alpha)\}$.

\item[c)] there exists a word
$\{\widetilde{t}^{t_j,t_{i+1}}\}=(T^{i_j,i_{j+1}}_{k(i_j,i_{j+1})},\cdots,T^{i_j,i_{j+1}}_1)$
with letters in $GL(R^d)$ $\varepsilon$ close to
$\{t^{i_j,t_{j+1}}\}$, let
$T^{i_j,i_{j+1}}=T^{i_j,i_j+1}_{k(i_j,i_j+1)}\cdot\cdots\cdot
T^{i_j,i_j+1}_{1}$, then \begin{equation*}
T^{i_j,i_{j+1}}(E^s(q_{i_{j}}))=E^s(q_{i_{j+1}}),\quad,T^{i_j,i_{j+1}}(E^u(q_{i_{j}}))=E^u(q_{i_{j+1}}).
\end{equation*}
\end{itemize}

We say $H(p,f)$ has transition property if $H(p,f)$ has
$\varepsilon$-transition property for any $\varepsilon>0$.

\begin{lem} \label{4.24}(\cite{1BDP}) $f\in C^1(M)$, suppose $p$ is an index $i$
($i\neq 0,d$) hyperbolic periodic point of $f$, then $H(p,f)$ has
transition property.\end{lem}

\begin{lem}\label{4.25} $f\in R$, suppose $p$ is an index $i$ ($i\neq
0,d$) hyperbolic periodic point of $f$ and $H(p,f)$ is not
trivial. Suppose there exists a family of periodic point
$\{p_n\}^\infty _{n=1}$ with index $i$ in $H(p,f)$ homoclinic
related with $p$ and $l_n\longrightarrow \infty$,
$\lambda_n\longrightarrow 1^-$ such that
$\pi_{l_n}(p_n)\longrightarrow \infty$ and $\prod
\limits_{j=0}^{\pi_{l_n}(p_n)-1} \Vert
Df^{l_n}|_{E^s(f^{jl_n}(p_n))}\Vert
\geq\lambda_n^{\pi_{l_n}(p_n)}$, then $H(p,f)$ is an index $i-1$
fundamental limit.\end{lem}

\Pf: We claim that we can find $q_n(g_n)$ is periodic point of
$g_n$ with index $i$ such that:
\begin{itemize}

\item[1)] $\lim \limits_{n\rightarrow \infty}g_n=f$.

\item[2)] $Orb_{g_n}(q_n)$ is periodic orbit of $f$ also
($f|_{Orb_{g_n}(q_n)}=g_n|_{Orb_{g_n}(q_n)}$), so we just denote
it $Orb(q_n)$, then we have $Orb(q_n)\subset H(p,f)$ and $\lim
\limits_{n\rightarrow \infty}Orb(q_n)=H(p,f)$.

\item[3)] $\lim \limits_{n\rightarrow \infty}
\frac{\pi(q_n)}{l_n}\longrightarrow \infty$

\item[4)]  $\prod \limits_{j=0}^{[\frac{\pi(q_n)}{l_n}]-1} \Vert
Dg_n^{l_n}|_{E^s_{g_n}(g_n^{jl_n}(q_n))} \Vert\geq
\lambda_n^{[\frac{\pi(q_n)}{l_n}]}$
\end{itemize}

\noindent{\bf Proof of the claim}: Choose
$\varepsilon_n\longrightarrow 0^+$, let's fix $n_0$ at first and
choose an $\varepsilon>0$ such that
$\lambda_{n_0}+2\varepsilon<1$. There exists $N_0\gg n_0$ such
that for any $n\geq N_0$, we'll have $l_n\gg l_{n_0}$ and
$\lambda_n>\lambda_{n_0}+2\varepsilon$, then by $\prod
\limits_{j=0}^{\pi_{l_n}(p_n)-1} \Vert
Df^{l_n}|_{E^s(f^{jl_n}p_n)}\Vert \geq
\lambda_n^{\pi_{l_n}(p_n)}$, we have $\prod
\limits_{j=0}^{ml_{n_0}\pi_{l_n}(p_n)-1} \Vert
Df^{l_n}|_{E^s(f^{jl_n}p_n)}\Vert \geq
\lambda_n^{ml_{n_0}\pi_{l_n}(p_n)}$ for $m\geq 1$, then we get
\begin{equation}\label{14}
\prod \limits_{j=0}^{ml_n\pi_{l_n}(p_n)-1} \Vert
Df^{l_{n_0}}|_{E^s(f^{jl_{n_0}}p_n)}\Vert \geq
(\lambda_{n_0}+2\varepsilon)^{ml_{n_0}\pi_{l_n}(p_n)} \;
\text{for}\  m\geq 1.\end{equation}

Since $f\in R$, there exists a family of periodic points
$\{q_i^\prime\}^N_{i=1}$ with index $i$, which are
$\varepsilon_{n_0}$-dense in $H(p,f)$ and they are homoclinic
related with $p$ and $\{p_n\}^\infty_{n=1}$. Now use
$\varepsilon_{n_0}$-transition property for
$\{q_0^\prime(=p_{N_0}),q_1^\prime,\cdots,q_N^\prime\}$, then for
$\{i,j\} \in\{0,1,\cdots,N\}^2$, there exists finite orbit
$[t^{i,j}]=(t^{i,j}_{k(i,j)},\cdots,t^{i,j}_{1})$ such that for
$l=(0,1,\cdots,N)$ and $\alpha_m=(m\cdot l_{n_0},1,\cdots,1)$, the
pseudo orbit
$[w(l,\alpha_m)]=[t^{N,0}]\cdot[q^\prime_N]\cdot\cdots\cdot[t^{0,1}]\cdot[q^\prime_0]^{m\cdot
l_{N_0}\frac{l_{N_0}\cdot \pi_{l_{N_0}}(p_{N_0})}{\pi(p_{N_0})}}$
is an $\varepsilon_{n_0}$-pseudo orbit and is
$\varepsilon_{n_0}$-shadowed by periodic orbit $[q(l,\alpha_m)]$
whose index is $i$, where $Orb(q(l,\alpha_m)))\subset H(p,f)$ and
$\{q(l,\alpha_m)\}$ is $\varepsilon_{n_0}$-near
$\{w(l,\alpha_m)\}$.

Consider the word
$\{\widetilde{w}(l,\alpha_m)\}=\{\widetilde{t}^{N,0}\}\cdot\{q^\prime_N\}\cdot\cdots\cdot\
\{\widetilde{t}^{N,0}\}\cdot \{q^\prime_0\}^{ml_0}$, it's
$\varepsilon_{n_0}$ near $\{w(l,\alpha_m)\}$, so
$\{\widetilde{w}(l,\alpha_m)\}$ is $2\varepsilon_{n_0}$ near with
$\{q(l,\alpha_m)\}$. Now use lemma \ref{Franks} (Franks lemma), we
can get a $C^1$ diffeomorphism $g_{(l,\alpha_m)}$ such that
$d(g_{(l,\alpha_m)},f)<2\varepsilon_{n_0}$, $Orb_f(q(l,\alpha_m))$
is also orbit of $g_{(l,\alpha_m)}$, and
$\{Dg_{(l,\alpha_m)}|_{Orb(q(l,\alpha_m))}\}=\{\widetilde{w}(l,\alpha_m)\}$.
By $c)$ of transition property, $E^{s(u)}_f(q^\prime_0)$ is
invariant bundle of $\{\widetilde{w}(l,\alpha_m)\}$, so they are
invariant bundle of $g_{l,\alpha_m}$, that means
$Dg^{\pi(q(l,\alpha_m))}_{(l,\alpha_m)}(E^s_f(q_0^\prime))=E^s_f(q_0^\prime)$
and
$Dg^{\pi(q(l,\alpha_m))}_{(l,\alpha_m)}(E^u_f(q_0^\prime))=E^u_f(q_0^\prime)$.
It's easy to know when $m$ is big enough, $E^{s(u)}_f(q_0^\prime)$
is stable(unstable) bundle for $g_{(l,\alpha_m)}$, so when $m$ is
big enough, $q_{(l,\alpha_m)}$ would be an index $i$ hyperbolic
periodic point of $g_{(l,\alpha_m)}$.

Now choose $m$ big enough and let $q_{n_0}=q(l,\alpha_m)$,
$g_{n_0}=g_{(l,\alpha_m)}$, it's easy to know $q_{n_0}$, $g_{n_0}$
satisfy $1)$, $2)$. About $3)$, let's notice that $\pi(q_n)\geq
ml_{n_0}$ and $m$ can be chosen arbitrary big. $4)$ comes from
\eqref{14} and $m$ is big enough. \qed

Now let's continue the proof of lemma \ref{4.25}, by the above
claim and lemma \ref{4.16}, for any $\varepsilon>0$ and $N>0$,
there exist an $n_0>N$ and $g_{n_0}^\prime$ is
$\varepsilon$-perturbation of $g_{n_0}$ such that $Orb(q_{n_0})$
is index $i-1$ periodic orbit of $g_{n_0}^\prime$ and
$Orb(q_{n_0})$ is $\varepsilon_{n_0}$-dense in $H(p,f)$. Since
$\varepsilon$ and $\varepsilon_{n_0}$ can be arbitrarily small, we
get that $\lim \limits_{n\rightarrow \infty}g_{n_j}^\prime=f$,
$Orb(q_{n_j})$ is index $i-1$ periodic orbit of $g_{n_j}^\prime$
and $\lim \limits_{j\rightarrow \infty}Orb(q_{n_j})=H(p,f)$, so
$H(p,f)$ is an index $i-1$ fundamental limit. \qed

  Then main result of this subsection is the following corollary.

\begin{cor} \label{4.26} $f\in R$, $C$ is a chain recurrent class of
$f$, $\Lambda$ is compact invariant set of $f$ with index
$i-(l,\lambda)$ dominated splitting $E^{cs}\oplus F^{cu}$ ($1\leq\
i\leq\ d$) and assume they satisfy all the assumption of weakly
selecting lemma, then $C$ contains index $i$ periodic point and
$C$ is an index $i-1$ funadamental limit.\end{cor}

\Pf: It's just a corollary from Lemma \ref{4.21} (weakly selecting
lemma) and lemma \ref{4.25}.    \qed

\subsection{Proof of lemma \ref{4.3}}

\Pf: When $\Lambda$ is trivial ($\# (\Lambda)<\infty$), $\Lambda$
is a periodic orbit, since $\Lambda$ is an index $j_0$-fundamental
limit, it should be an index $j_0$ hyperbolic periodic orbit, so
$C$ contains an index $j_0$ periodic point and it's an index $j_0$
fundamental limit.

Now we suppose $\Lambda$ is not trivial, by generic property $5$
of proposition 3.1, there exists a family of index $j_0$ periodic
points $\{p_n(f)\}^\infty_{n=1}$ such that $\lim
\limits_{n\rightarrow \infty}Orb(p_n(f))=\Lambda$. Since $\Lambda$
is not trivial, we have $\pi(p_n(f))\longrightarrow \infty$.

If $\Lambda$ isn't an index $j_0+1$ fundamental limit, we know
that $\{p_n(f)\}$ is index $j_0$ stable, then by lemma \ref{4.6}
(Gan's lemma), there exits a subsequence
$\{p_{n_i}(f)\}^\infty_{i=1}$ such that $p_{n_i}(f)$ and
$p_{n_j}(f)$ are homoclinic related, so
$H(p_{n_1},f)=H(p_{n_2},f)=\cdots$, especially, by $\lim
\limits_{n\rightarrow \infty}Orb(p_n(f))=\Lambda$, we know that
$\Lambda \subset H(p_{n_1},f)$, by generic property $6)$ of
proposition 3.1, $C=H(p_{n_1},f)$, so $C$ contains index $j_0$
periodic point and it's an index $j_0$ fundamental limit.

So from now, we suppose $\Lambda$ is an index $j_0+1$ fundamental
limit also, then $\Lambda \subset P^*_{j_0}\bigcap P^*_{j_0+1}$,
since $f$ is far away from tangency, by proposition 2.1, $\Lambda$
has an index $j_0$ dominated splitting
$E^{cs}_{j_0}(\Lambda)\oplus E^{cu}_{j_0+1}(\Lambda)$ and an index
$j_0+1$ dominated splitting $E^{cs}_{j_0+1}(\Lambda)\oplus
E^{cu}_{j_0+2}(\Lambda)$. Let
$E_1^c(\Lambda)=E^{cu}_{j_0+1}(\Lambda) \bigcap
E^{cs}_{j_0+1}(\Lambda)$, then on $\Lambda$ we have the following
dominated splitting: $T|_{\Lambda}M=E^{cs}_{j_0}(\Lambda)\oplus
E^{c}_{1}(\Lambda) \oplus E^{cu}_{j_0+2}(\Lambda)$. Since
$C\bigcap P^*_j= \phi$ for $j<j_0$, by proposition 2.2,
$E^{cs}_{j_0}$ is in fact contracting, so we prefer denoting it
$E^s_{j_0}$. Now on $\Lambda$ we have the dominated splitting
$T|_{\Lambda}M=E^{s}_{j_0}(\Lambda)\oplus E^{c}_{1}(\Lambda)
\oplus E^{cu}_{j_0+2}(\Lambda)$.

\begin{rk}\label{4.27} Since $\Lambda$ is index $j_0$ fundamental limit,
$E^c_1(\Lambda)$ is not contracting, that means that the bundle
$(E^s_{j_0}\oplus E^c_1)|_\Lambda$ is not contracting
also.\end{rk}

When $j_0+1=d$, especially, the dominated splitting on $\Lambda$
should be $T|_\Lambda M=E^s_{j_0}(\Lambda)\oplus E^c_1(\Lambda)$.
In this case, if $\Lambda$ is not minimal, there exists an $x_0\in
\Lambda$ such that $\omega(x_0)\varsubsetneq \Lambda$. By the
definition of $\Lambda$ and $j_0=d-1$, $\omega(x_0)$ is an index
$d$ fundamental limit but not index $j$ fundamental limit for
$j<d$. With the generic property $(5)$ of proposition 3.1,
$\omega(x_0)$ can be converged by a family of sinks $\{p_n(f)\}$,
by remark 4.4, $\pi(p_n(f))$ should be bounded ( If it's not
bounded, there exist $p_{n_0}(f)$ and
$g_{n_0}\overset{C^1}{\sim}f$ such that
$g_{n_0}|_{Orb_f(p_{n_0}(f))}=f|_{Orb_f(p_{n_0}(f))}$ and
$Orb(p_{n_0}(f))$ is a periodic orbit of $g$ with index smaller
than $d$, that means $\omega(x_0)$ is an fundamental limit with
index smaller than $d$, it's a contradiction). That means
$\omega(x_0)$ is trivial, so it's a periodic orbit. Since $f$ is a
Kupuka-Smale diffeomorphism and $\omega(x_0)$ is an index $d$
fundamental limit, we can know that $\omega(x_0)$ is an index $d$
hyperbolic periodic orbit, then $C$ contains a sink, it means $C$
itself is just the orbit of sink and $C=\omega(x_0)$, that's a
contradiction with $C$ is not trivial, so we proved $\Lambda$ is
minimal when $j_0+1=d$.

Now we just consider $j_0+1<d$, we claim that with all the
assumptions above on $\Lambda$, then either $\Lambda$ is minimal,
or $C$ contains periodic points with index $j_0+1$ and $C$ is an
index $j_0$ fundamental limit.\\

\noindent{\bf Proof of claim}: Suppose $\Lambda$ is not minimal,
it means that there exists $x_0\in \Lambda$ such that
$\omega(x_0)\neq \Lambda$. Consider the set of compact chain
recurrent subset of $\Lambda$: $\{\Lambda_\alpha:\ \Lambda_\alpha
\varsubsetneq \Lambda\}_{\alpha\in \cal A}$, since $\omega(x_0)\in
\{\Lambda_\alpha\}_{\alpha\in \cal A}$, $\cal A \neq \phi$, by
generic property $(4)$ of proposition 3.1, $\Lambda_\alpha$ is a
fundamental limit. By the definition of $j_0$ and $\Lambda$,
$\Lambda_\alpha$ is an index $j_\alpha$ fundamental limit with
$j_\alpha \geq j_0+1$. Denote $\cal B$=$\{\beta\in \cal A$,
$\Lambda_\beta$ is not an index $j$ fundamental limit for
$j>j_0+1\}$.

\begin{rk}\label{4.28}: For any $\beta \in \cal B$, $\Lambda_\beta$ is an
index $j_0+1$ fundamental limit, on $\Lambda_\beta$ we have an
index $j_0+1$ dominated splitting
$E^{cs}_{j_0+1}(\Lambda_\beta)\oplus
E^{cu}_{j_0+2}(\Lambda_\beta)$. Since we have $\Lambda_\beta
\bigcap P^*_j \neq \phi$ for all $j \neq j_0+1$, by proposition
2.2, the index $j_0+1$ dominated splitting is in fact a hyperbolic
splitting, that means $\Lambda_\beta$ is a hyperbolic set.\end{rk}

  Now we divide the proof of the claim to three subcases: $\#(B)=0$, $\#(B)=N_1<\infty$ and
  $\#(B)=\infty$.\\

\noindent{Case A}: $\#(B)=0$.\\

  That means for all $\alpha\in \cal A$, $\Lambda_\alpha$ is an index $j_\alpha$ fundamental limit for some
  $j_\alpha>j_0+1$.

  Now we need the following two results.

\begin{lem}\label{4.29} (\cite{1W4}) Assume $f\in R$, let $\Lambda$ be an index
$i$ fundamental limit of $f$ ($1\leq\ i\leq \ d-1$),
$E^{cs}_i(\Lambda)\oplus E^{cu}_{i+1}(\Lambda)$ is an index
$i-(l,\lambda)$ dominated splitting on $\Lambda$ given by
proposition 2.1, then
\begin{itemize}
       \item[1)] either for any $\mu \in(\lambda,1)$, there exists $c\in \Lambda$ such that
       $\prod \limits_{j=0}^{n-1} \Vert Df^l|_{E^{cs}_i(f^{jl}c)}\Vert \leq \mu^n$ for
       $n\geq 1$,
       \item[2)] or $E^{cs}_i$ splits into a dominated splitting $V^{cs}_{i-1}\oplus V^{c}_{1}$ with $dim(V^c_1)=1$ such
       that for any $\mu \in(\lambda,1)$, there is $c^\prime\in \Lambda$ such that
       $\prod \limits_{j=0}^{n-1} \Vert Df^l|_{V^{cs}_{i-1}(f^{jl}c^\prime)}\Vert \leq \mu^n$ for all $n \geq 1$.
\end{itemize}
\end{lem}

\begin{lem} \label{4.30} Let $\Lambda$ be an invariant compact set of
$f$, with two dominated splitting $E^{cs}\oplus F^{cu}$ and
$\widetilde{E}^{cs}\oplus\widetilde{F}^{cu}$, if $dim(E^{cs})\leq
dim(\widetilde{E}^{cs})$, then $E^{cs}\subset
\widetilde{E}^{cs}$.\end{lem}

Choose $\mu_0\in(\lambda,1)$, since $\Lambda_\alpha$ is an index
$j_\alpha$ fundamental limit, proposition 2.1 gives an index
$j_\alpha-(l,\lambda)$ dominated splitting
$E^{cs}_{j_\alpha}\oplus F^{cu}_{j_\alpha+1}$ on $\Lambda_\alpha$.

If 1) of lemma \ref{4.29} is true for $\Lambda_\alpha$, then there
exists $c\in\Lambda_\alpha$ such that $\prod
\limits_{j=0}^{n-1}\Vert Df^l|_{E^{cs}_{j_\alpha}(f^{jl}c)}\Vert
\leq\mu^n_0$ for $n\geq 1$. On $\Lambda_\alpha$ we have another
dominated splitting $(E^{s}_{j_0}\oplus E^{c}_{1})\oplus
E_{cu}^{j_0+2}$ induced from $\Lambda$. Since
$dim(E^{s}_{j_0}\oplus
E^{c}_{1})=j_0+1<j_\alpha=dim(E^{cs}_{j_\alpha})$, be lemma
\ref{4.30}, $E^{s}_{j_0}\oplus E^{c}_{1}\subset
E^{cs}_{j_\alpha}$, so we have $\prod \limits_{j=0}^{n-1}\Vert
Df^l|_{E^{s}_{j_0}\oplus E^{c}_{1}(f^{jl}c)}\Vert \leq\mu_0^n$ for
$n\geq 1$.

  If 2) of lemma \ref{4.29} is true for $\Lambda_\alpha$, then there exists
$c^\prime$ such that $\prod \limits_{j=0}^{n-1} \Vert
Df^l|_{V^{cs}_{j_\alpha-1}(f^{jl}c^\prime)}\Vert \leq \mu^n_0$ for
$n\geq 1$, recall that $dim(E^{s}_{j_0}\oplus E^{c}_{1})=j_0+1\leq
j_\alpha-1=dim(V^{cs}_{j_\alpha-1})$, by lemma \ref{4.30},
$E^{s}_{j_0}\oplus E^{c}_{1}\subset
V^{cs}_{j_\alpha}(\Lambda_\alpha)$, so we have $\prod
\limits_{j=0}^{n-1}\Vert Df^l|_{E^{s}_{j_0}\oplus
E^{c}_{1}(f^{jl}c^\prime)}\Vert \leq\mu_0^n$ for $n\geq 1$.

\begin{rk} \label{4.31}: By the above arguments, we know that for any
$\alpha \in \cal A \setminus B$, and for any
$\mu_0\in(\lambda,1)$, there exists $c\in \Lambda_\alpha$ such
that \begin{equation}\label{a15}
 \prod
\limits_{j=0}^{n-1}\Vert Df^l|_{E^{s}_{j_0}\oplus
E^{c}_{1}(f^{jl}c)}\Vert \leq\mu_0^n\quad \text{for}\  n\geq 1.
\end{equation}\end{rk}

By remark \ref{4.27}and remark \ref{4.31}, the index
$j_0+1-(l,\lambda)$ dominated splitting $(E^s_{j_0}\oplus
E^c_1)\oplus E^{cu}_{j_0+2}$ on $\Lambda$ satisfies all the
conditions of weakly selecting lemma, by corollary \ref{4.26}, $C$
contains index $j_0+1$ periodic
point and $C$ is an index $j_0$ fundamental limit.\\

\noindent{Case B}: $\#(B)=N_1<\infty$

Let $\cal B$=$\{\beta_1,\cdots,\beta_{N_1}\}$, fix
$\lambda<\mu_0<1$, then by the argument in case A, for any $\beta
\in \cal A\setminus B$, there exists $c\in \Lambda$ satisfies
\eqref{a15}.

 For $\beta_i\in \cal B$, $\Lambda_{\beta_i}$ should be a hyperbolic set where the
bundle $E^s_{j_1}\oplus E^c_1|_{\Lambda_{\beta_i}}$ is a
contracting bundle, so there exists $l^\prime$ such that for any
$x\in \Lambda_{\beta_i}$, $\Vert Df^{l^\prime}|_{(E^s_{i_0}\oplus
E^c_1)(x)}\Vert <1/2$.

 Let $l_0=l\cdot\l^\prime$ and $1>\mu_1>\max \{\mu_0,\frac{1}{2}\}$, then for any
$\Lambda_\alpha$ ($\alpha\in \cal A$), there exists a point $c\in
\Lambda_{\alpha}$ such that $\prod \limits_{j=0}^{n-1}\Vert
Df^{l_0}|_{E^{s}_{j_0}\oplus E^{c}_{1}(f^{jl_0}c)}\Vert
\leq\mu_1^n$. With remark \ref{4.27}, the index $j_0+1-(l,
\lambda)$ dominated splitting $(E^s_{j_0}\oplus E^c_1)\oplus
E^{cu}_{j_0+2}$ on $\Lambda$ satisfies all the conditions of
weakly selecting lemma, by corollary \ref{4.26}, $C$ contains
index $j_0+1$ periodic point and
$C$ is an index $j_0$ fundamental limit.\\

\noindent{Case C}: $\#(B)=\infty$

  In remark \ref{4.28}, we have shown that for any $\beta\in \cal B$, $\Lambda_\beta$ is a hyperbolic chain recurrent set with
index $j_0+1$. Then there exists a family of periodic points
$\{p_{\beta,n}\}^\infty_{n=1}$ in $C$ with index $j_0+1$ and $\lim
\limits_{n\rightarrow \infty}Orb(p_{\beta,n})=\Lambda_\beta$ (by
shadowing lemma). If $\Lambda_\beta$ is trivial, that means it's
an index $j_0+1$ periodic orbit, we can let
$Orb(p_{\beta,n})=\Lambda_\beta$ for $n\geq 1$; if $\Lambda_\beta$
is not trivial, we can let $\pi(p_{\beta,n})\longrightarrow
\infty$.

  We have the following two subcases.\\
\begin{itemize}
  \item[$\bullet$] Subcase $C.1$: There exists $\delta>0$ such that for any $\Lambda_\beta,\ \beta\in \cal B$,
  there exists a family of periodic points $\{p_{\beta,n}\}^\infty_{n=1}$ such that
  $\lim \limits_{n\rightarrow \infty}Orb(p_{\beta,n})=\Lambda_\beta$ and
  $\mid Df^{\pi(p_{\beta,n})}|_{E^c_1(p_{\beta,n})}\mid <e^{-\delta\pi(p_{\beta,n})} $.

  \item[$\bullet$] Subcase $C.2$: For any $\frac{1}{m}>0$, there exist
$\beta_m \in \cal B$ and a family of periodic points
$\{p_{\beta_m,n}\}_{n=1}^\infty$ satisfying $\lim
\limits_{n\rightarrow \infty}Orb(p_{\beta_m,n})=\Lambda_\beta$ and
$\mid Df^{\pi(p_{\beta_m,n})}|_{E^c_1(p_{\beta_m,n})}\mid
>e^{-\frac{1}{m}\pi(p_{\beta_m,n})} $.
\end{itemize}

  In the subcase $C.1$, let's fix $1>\mu_1>\mu_0>e^{-\delta}$. For $\beta\in \cal B$, recall that $dim(E^c_1(\Lambda))=1$
  and $\mid Df^{\pi(p_{\beta,n})}|_{E^c_1(p_{\beta,n})}\mid <e^{-\delta \pi(p_{\beta,n})}$,
  we'll get $\prod \limits_{i=0}^{\pi(p_{\beta,n})-1}\mid Df|_{E^c_1(p_{\beta,n})}\mid<e^{-\delta \pi(p_{\beta,n})}$,
  that means for any $s\geq 1$, we have
  $\prod \limits_{i=0}^{s\pi(p_{\beta,n})-1}\mid Df|_{E^c_1(p_{\beta,n})}\mid<e^{-s\delta \pi(p_{\beta,n})}$ for
$s\geq 1$. By lemma \ref{4.10} (Pliss lemma) there exists
$x_{\beta,n}\in Orb(p_{\beta,n})$ such that $\mid
Df^s|_{E^c_1(x_{\beta,n})}\mid=\prod \limits_{i=0}^{s-1}\mid
Df|_{E^c_1(f^i(x_{\beta,n}))}\mid<\mu_0^s$ for $s \geq 1$. Suppose
$\lim \limits_{n\rightarrow \infty}x_{\beta,n}\longrightarrow
c_\beta$ where $c_\beta \in \Lambda_\beta$, then $\prod
\limits_{i=0}^{s-1}\mid Df|_{E^c_1(f^i(c_\beta))}\mid<\mu_0^s$ for
$s\geq 1$. Notice that $E^s_{j_0}|_\Lambda$ is dominated by
$E^c_1|_\Lambda$ and $\mu_1>\mu_0$, there exists $l^\prime\gg 1$
doesn't depend on $\beta$ such that $\prod
\limits_{i=0}^{t-1}\Vert Df^{l^\prime}|_{E^c_1\oplus
E^s_{j_0}(f^{il^\prime}(c_\beta))}\Vert<\mu_1^t$ for $t \geq 1$.

    For $\alpha \in \cal A\setminus B$, by the argument in case A, there exists
$c_\alpha \in \cal A_\alpha$ such that $\prod
\limits_{i=0}^{t-1}\Vert Df^{l_0}|_{E^c_1\oplus
E^s_{j_0}(f^{il_0}(c_\alpha))}\Vert<\mu_1^t$ for $t\geq 1$.

    Let $l_1=l^\prime\cdot l_0$, then for any $\alpha \in \cal A$, there exists $c_\alpha \in \cal A$
such that $\prod \limits_{i=0}^{t-1}\Vert Df^{l_1}|_{E^c_1\oplus
E^s_{j_0}(f^{il_1}(c_\alpha))}\Vert<\mu_1^t$ for $t\geq 1$. With
remark \ref{4.27}, the index $j_0+1-(l,\lambda)$ dominated
splitting $(E^s_{j_0}\oplus E^c_1)\oplus E^u_{j_0+2}$ on $\Lambda$
satisfies all the conditions of weakly selecting lemma. By
Corollary \ref{4.26}, $C$ contains index $j_0+1$ periodic point
and $C$ is an index $j_0$ fundamental limit.

  In the subcase $C.2$, since $\Lambda_{\beta_m}$ is a hyperbolic set, we can always suppose
  $\{p_{\beta_m,n}\}^\infty_{n=1}$ is homoclinic
related with each other and $p_{\beta_m,n}\in C$, so $C$ contains
index $j_0+1$ periodic points. Now we'll show $C$ is an index
$j_0$ fundamental limit also.

  We claim that there exists a subsequence $\{\beta_{m_t}\}^\infty_{t=1}\subset \{\beta_m\}$ and for every
$\beta_{m_t}$ there exists $p_{\beta_{m_t},n_t}\in
\{p_{\beta_{m_t},n}\}^\infty_{n=1}$ such that $\lim
\limits_{t\rightarrow
\infty}\pi(p_{\beta_{m_t},n_t})\longrightarrow \infty$.\\

\noindent{\bf Proof of the claim}: Let ${\cal B}_0$=$\{ \beta_m:\
\Lambda_{\beta_m}$ is given in subcase $C.2$ and
$\Lambda_{\beta_m}$ is not trivial. $\}$

   If $\#({\cal B}_0)=\infty$, then for any $\beta_{m_t}\in \cal B$$_0$,
by $\Lambda_{\beta_{m_t}}$ is not trivial, we'll have $\lim
\limits_{n\rightarrow \infty}\pi(p_{\beta_{m_t},n})\longrightarrow
\infty$, so when n is big enough, we can let
$\pi(p_{\beta_{m_t},n})$ arbitrarily big.

   If $\#({\cal B}_0)<\infty$, then for $\beta_m \notin \cal B$$_0$, $\Lambda_{\beta_m}$
is an index $j_0+1$ periodic orbit and $Orb(p_{\beta_m,n})\equiv
\Lambda_{\beta_m}$ for $n\geq 1$. Since $f$ is a Kupka-Smale
diffeomorphism, the number of periodic points with fixed boundary
of period should be finite, by the fact $\# ({\cal B\setminus \cal
B}_0)=\infty$, there are infinite of $m$ such that $\Lambda_m$ is
index $j_0+1$ periodic orbits, then we can choose
$\Lambda_{\beta_m}$ is an index $j_0+1$ periodic orbit with
arbitrarily big period.  \qed

Now for simiplicity, we denote $p_{\beta_{m_t},n_t}$ by
$p_{\beta_m,n_m}$.

For $\{p_{\beta_m,n_m}\}^\infty_{m=1}$, we have $\lim
\limits_{m\rightarrow \infty}\pi(p_{\beta_m.n_m})\longrightarrow
\infty$ and
\begin{equation}\label{15}
|Df^{\pi(p_{\beta_m,n_m})}|_{E^c_1(p_{\beta_m,n_m})}|>e^{-
\frac{1}{m}\pi(p_{\beta_m,n_m})}.
\end{equation}Choose
$\{l_m\}^\infty_{m=1}$ carefully, we'll have $\lim
\limits_{m\rightarrow \infty}l_m\longrightarrow \infty$, $\lim
\limits_{m\rightarrow \infty}
\frac{\pi(p_{\beta_m,n_m})}{l_m}\longrightarrow \infty$ and
$\frac{l_m}{m}\longrightarrow 0^+$ (after replacing
$\{p_{\beta_m,n_m}\}^\infty_{m=1}$ by a subsequence, we can always
do this). Since $\pi_{l_m}(p_{\beta_m,n_m})\geq
\frac{\pi(p_{\beta_m,n_m})}{l_m}$, we have
\begin{equation}\label{16} \lim \limits_{m\rightarrow
\infty}\pi_{l_m}(p_{\beta_m,n_m})\longrightarrow \infty.
\end{equation}

  By \eqref{15} and the fact $l\cdot \pi_l(p)$ is always a multiple of $\pi(p)$ for any period point $p$
and $l\geq 1$, we have
$$|Df^{l_m\cdot\pi_{l_m}(p_{\beta_m,n_m})}|_{E^c_1(p_{\beta_m,n_m})}|>e^{-\frac{1}{m}l_m\cdot\pi_{l_m}(p_{\beta_m,n_m})},$$
it's equivalent with $$\prod
\limits_{i=0}^{\pi_{l_m}(p_{\beta_m,n_m})-1}\Vert
Df^{l_m}|_{E^c_1(f^{il_m}(p_{\beta_m,n_m}))}\Vert \geq
e^{-\frac{l_m}{m}\cdot \pi_{l_m}(p_{\beta_m,n_m})},$$ then we get
$$\prod \limits_{i=0}^{\pi_{l_m}(p_{\beta_m,n_m})-1}\Vert
Df^{l_m}|_{(E^c_1\oplus
E^s_{j_0})(f^{il_m}(p_{\beta_m,n_m}))}\Vert \geq
e^{-\frac{l_m}{m}\cdot \pi_{l_m}(p_{\beta_m,n_m})},$$ since $\lim
\limits_{m\rightarrow \infty}\frac{l_m}{m}\longrightarrow 0^+$ and
by \eqref{16}, lemma \ref{4.25} tells us $C$ is an index $j_0$
fundamental limit, this finishes the proof of the claim. \qed

Now let's continue the proof of lemma \ref{4.3}, by the above
argument, we can suppose $\Lambda$ is minimal, not trivial, it's
an index $j_0$ and $j_0+1$ fundamental limit with dominated
splitting $E^s_{j_0}\oplus E^c_1\oplus E^{cu}_{j_0+2}|_\Lambda$ where
$E^{cu}_{j_0+2}(\Lambda)\neq \phi$.

  If $E^{cu}_{j_0+2}(\Lambda)$ is not expanding, by lemma \ref{4.22},
we can know that there exists a point $b\in \Lambda$ such that
$\prod \limits_{i=0}^{n-1}\Vert
Df^{-l}|_{E^{cu}_{j_0+2}(f^{(i+1)l}b)}\Vert \geq 1$, since
$(E^s_{j_0}\oplus E^c_1)\oplus E^{cu}_{j_0+2}|_\Lambda$ is an
index $j_0+1-(l,\lambda)$ fundamental limit, it means that $$\prod
\limits_{i=0}^{n-1}\Vert Df^l|_{E^s_{j_0}\oplus
E^c_1(f^{il}(b))}\Vert \cdot\prod \limits_{i=0}^{n-1}\Vert
Df^{-l}|_{E^{cu}_{j_0+2}(f^{(i+1)l}(b))}\Vert \leq \lambda^n,\quad
\text{for}\ n\geq 1,$$ so $\prod \limits_{i=0}^{n-1}\Vert
Df^l|_{E^s_{j_0}\oplus E^c_1(f^{il}(b))}\Vert\leq \lambda^n$ for
all $n \geq 1$. Since $\Lambda$ is minimal, the index $j_0+1$
dominated splitting on $\Lambda$ satisfies strong tilda condition,
by remark \ref{4.27}, it also satisfies the non-hyperbolic
condition, so it satisfies all the conditions of weakly selecting
lemma, then by corollary \ref{4.26}, $C$ contains index $j_0+1$
periodic point and it's an index $j_0$ fundamental limit. \qed

\section{Proof of theorem 1}
  In order to prove theorem 1, we need the following lemma whose proof has been postponed to the end of this section.

\begin{lem}\label{5.1} Let $f\in R$, $C$ is any non-trivial chain recurrent class
of $f$, suppose $\Lambda\subset C$ is a non-trivial minimal set
with a codimension-1 partial hyperbolic splitting $T_\Lambda
M=E^c_1\oplus E^u_2$ where $dim(E^c_1|_\Lambda)=1$ and
$E^c_1(\Lambda)$ is not contracting, then $C$ is a homoclinic
class containing index 1 periodic point and $C$ is an index 0
fundamental limit.
\end{lem}

\begin{rk}\label{5.2} in \cite{1BGW}, they show that for $f\in R$, if $C$ is a chain
recurrent class of $f$ with a codimension-1 dominated splitting
$T_CM=E^c_1\oplus E^u_2$ where $dim(E^c_1|_C)=1$ and $E^c_1|_C$ is
not hyperbolic, then $C$ should be a homoclinic class. We
generalize this result to minimal set with Crovisier's work on
central curves.
\end{rk}

\noindent{\bf Proof of theorem 1}: Suppose $C\bigcap P^*_0\neq
\phi$, let $\Lambda$ be an minimal index 0 fundamental limit, then
$\Lambda$ is not trivial ( if $\Lambda$ is trivial, $\Lambda$
should be an orbit of source, then $C$ itself is source also, that
contradicts with $C$ is not trivial)). By lemma \ref{4.3}, either
$C$ is a homoclinic class containing index 1 periodic point and
$C$ is an index 0 fundamental limit or $\Lambda$ is a non-trivial
minimal set with codimension-1 partial hyperbolic splitting
$T_\Lambda M=E^c_1\oplus E^u_2$ where $E^c_1|_\Lambda$ is not
trivial. In the first case we've proved theorem 1, in the second
case, by lemma \ref{5.1}, we also prove theorem 1.\qed

  In $\S$5.1, we'll introduce some properties for codimension-1
partial hyperbolic splitting set, in $\S$5.2 we'll introduce
Crovisier's central model for the invariant compact set with
partial hyperbolic splitting whose central bundle is 1-dimension
and non-hyperbolic. In $\S$5.3 I'll give the proof of lemma
\ref{5.1}.
\subsection{Some properties for codimension-1 partial hyperbolic splitting}

Let $f\in R$, $\Lambda$ is a given non-trivial minimal set of $f$
with a codimension-1 partial hyperbolic splitting $T_\Lambda
M=E^u\oplus E^c_1$, where $dim(E^c_1(\Lambda))=1$ and the bundle
$E^c_1|_\Lambda$ is not hyperbolic. In this section we always
suppose the dominated splitting is 1-step and the bundle $E^u$ is
1-step expanding, it means that there exists $0<\lambda<1$ such
that for any $v^u\in E^u(x)$, $v^c\in E^c_1(x)$ where
$|v^u|=|v^c|=1$, $x\in \Lambda$, we have
$\frac{|Df(v^c)|}{|Df(v^u)|}<\lambda$, $|Df(v^u)|>\lambda^{-1}$. Fix a
small neighborhood $U_0$ of $\Lambda$, then the maximal invariant
set $\Lambda_0=\bigcap
\limits_{j=-\infty}^{\infty}f^j(\overline{U_0})$ has also a
codimension-1 partial hyperbolic splitting
$\widetilde{E^u}\oplus\widetilde{E^c_1}$, the dominated
splitting is 1-step and the bundle $\widetilde{E^u}|_{\Lambda_0}$
is also 1-step expanding. We say $E^c_1(\Lambda)$ has an
$f$-orientation if $E^c_1|_\Lambda$ is orientable and $Df$
preserves the orientation. If $E^c_1|_\Lambda$ has an
$f$-orientation, we choose $U_0$ small enough such that
$\widetilde{E^c_1}(\Lambda)$ has an $f$-orientation also.

Here we should notice the reader that in this section, all the
argument will take place just in $U_0$, and we can suppose $U_0$
is small enough such that it satisfies all the properties which we
need.

When $U_0$ is small enough, we can extend the bundle
$\widetilde{E^u}|_{\Lambda_0}$ and
$\widetilde{E^c_1}|_{\Lambda_0}$ to $\overline{U_0}$ such that for any $x\in
\overline{U_0}$, $T_xM=\widetilde{E^u}(x)\oplus\widetilde{E^c_1}(x)$, and if
$E^c_1|_\Lambda$ is orientable,
$\widetilde{E^c_1}|_{\overline{U}_0}$ is orientable also. In fact,
no matter $\widetilde{E^c_1}|_{\overline{U}_0}$ is orientable or
not, we can always locally define an orientation of
$\widetilde{E^c_1}|_{\overline{U}_0}$, it means that there exists
$\delta_0>0$ such that for any $x\in \overline{U_0}$, we can
give an orientation for the bundle
$\widetilde{E^c_1}|_{B_{\delta_0}(x)}\bigcap \overline{U}_0$.

For every point $x\in \overline{U_0}$, we define two kinds of
cones on its tangent space $C^i_a(x)=\{v|v\in T_xM,$ there exists
$v^\prime \in \widetilde{E^i}(x)$ such that
$d(\frac{v}{|v|},\frac{v^\prime}{|v^\prime|})<a\}_{i=c,u}$. When
$a$ small enough, $C^c_a\bigcap C^u_a=\phi$, $Df(C^u_a(x))\subset
C^u_a(f(x))$ and $Df^{-1}(C^c_a(x))\subset C^c_a(f^{-1}(x))$ for
$x\in \Lambda_0$.

We say a submanifold $D^i$ ($i=c,u$) tangents with cone $C^i_a$ if
$dim (D^i)=d-1$ when $i=u$ and $dim(D^i)=1$ when $i=c$ and for $x\in
D^i$, $T_xD^i\subset C^i_a(x)$. For simplicity, sometimes we call
it $i$-disk, especially when $i=c$, we just call $D^c$ a central
curve. We say an $i$-disk $D^i$ has centrer $x$ with size $\delta$
if $x\in D^i$, and respecting the Riemannian metric restricting on
$D^i$, the ball centered on $x$ with radius $\delta$ is in $D^i$.
We say an $i$-disk $D^i$ has center $x$ with radius $\delta$ if
$x\in D^i$, and respecting the Riemannian metric restricting on
$D^i$, the distance between any point $y\in D^i$ and $x$ is
smaller than $\delta$.

The following lemma shows some well-known results, it depends on a
simple fact: locally the splitting $\widetilde{E^c_1}\oplus
\widetilde{E^u}|_{\overline{U}_0}$ looks like linear. \cite{1BGW}
's subsection 4.1 gives many details about such view, from lemma
4.8 in \cite{1BGW}, it would be very easy to get the following
properties, so here we 'll not give a proof.

\begin{lem} \label{5.3}: Let $f\in R$, $\Lambda$ is a non-trivial minimal set of $f$ with
a codimension-1 partial hyperbolic splitting $T_\Lambda
M=E^c_1\oplus E^u$ where the bundle $E^c_1|_\Lambda$ is not
hyperbolic. $U_0,\delta_0,C^u_a,C^c_a$ are defined by the above
argument. Let $U$ be any small neighborhood of $\Lambda$
satisfying $\overline{U}\subset U_0$, there exist two
neighborhoods $U_2,U_1$ of $\Lambda$ such that $\Lambda \subset
U_2\subset \overline{U_2}\subset U_1 \subset \overline{U_1}\subset
U\subset U_0$ and there exist $a_0$ small enough and
$0<\delta_{1,3}<\delta_{1,2}<\delta_{1,1}<\delta_0/2$ such that
they satisfy the following properties:
\begin{itemize}
\item[P1£©]For any $x\in \overline{U_2}$, we have
$B_{2\delta_{1,1}}(x)\subset U_1$, and for any $x\in
\overline{U_1}$, we have $B_{2\delta_{1,1}}(x)\subset U$, then for
any $i$-disk $D^i$ ($i=c,u$) with center $x \in \overline{U_1}$
and radius $2\delta_{1,1}$ we'll have $D^i\subset U$.

\; For any $x\in \overline{U_1}$,
$\widetilde{E^c_1}|_{B_{2\delta_{1,1}}(x)}$ is orientable, we can
choose an orientation and call the direction right, then the
orientation of $\overline{E^c_1}|_{B_{2\delta_{1,1}}(x)}$ will
give an orientation for central curves in $B_{2\delta_{1,1}}(x)$.
We suppose $\delta_{1,1}$ is small enough such that any central
curve in $B_{2\delta_{1,1}}(x)$ will not intersect with itself.

\; For two points $y_1,y_2\in B_{2\delta_{1,1}}(x)$, we say $y_1$
is on the $x$-right of $y_2$ if there exists a central curve
$l\subset B_{2\delta_{1,1}}(x)$ connects $y_1$ and $y_2$ and in
$l$, $y_1$ is on the right of $y_2$. Then since any central curve
in $B_{2\delta_{1,1}}(x)$ is not self-intersection, $y_2$ is not
on $x$-right of $y_1$ anymore. Usually, we just simply call $y_1$
is on the right of $y_2$.

\item[P2)]Let $\Lambda_1=\bigcap \limits_{i=-\infty}^{\infty}f^i(\overline{U_1})$,
apply lemma \ref{4.13} on $\Lambda_1$, we can get the following
two kinds of submanifolds families: the local unstable manifolds
$W^{uu}_{loc}(x)\ _{x\in \Lambda_1}$ and the local central curves
$W^c_{loc}(x)\ _{x\in \Lambda_1}$.

\; Choose $\delta_{1,1}$ properly ( small enough) we can suppose
$W^i_{loc}(x)\ _{(i=uu,c)}$ has size $\delta_{1,1}$, let
$W^i_{\delta_{1,1}}(x)$ be the ball in $W^i_{loc}(x)$ with central
$x$ and radius $\delta_{1,1}$, then we have
$W^i_{\delta_{1,1}}(x)\ _{(x\in \Lambda_1,i=c,uu)}$ always
tangents with cone $C^i_{a_0}$.

\; In fact, for $\Lambda^+_1=\bigcap
\limits_{i=0}^{\infty}f^i(\overline{U_1})$, any $x\in \Lambda^+_1$
will have uniform size of unstable manifold
$W^{uu}_{\delta_{1,1}}(x)$ which tangents with cone
$C^{uu}_{a_0}$.

\item[P3)] By the property of strong unstable manifolds, for
$y_1,y_2\in\Lambda_1^+$, if we have
$W^{uu}_{\delta_{1,1}/2}(y_1)\bigcap W^{uu}_{\delta_{1,1}/2}(y_2)$
$\neq \phi$, then $y_1\in W^{uu}_{\delta_{1,1}}(y_2)$ and $y_2\in
W^{uu}_{\delta_{1,1}}(y_1)$. There exists $0<\lambda<1$ such that
for any smooth curve $l\subset W^{uu}_{\delta_{1,1}}(x)$ where
$x\in \Lambda^+_1$, we'll have $length(f^{-1}(l))<\lambda \cdot
length(l)$.

\item[P4)]For any central curve $D^c$ and $u$-disk $D^u$ in $U$ with
centers in $\Lambda_1$ and radius smaller than $2\delta_{1,1}$,
we have $\#\{z|\ z\in D^c\bigcap D^u\}\leq 1$. If $D^c\bigcap D^u
\neq \phi$, then they are transverse intersect with each other.

\item[P5)] For any $x\in \overline{U_1}$, $y\in
B_{\delta_{1,3}}(x)\bigcap \Lambda_1$, $D^i_{\delta_{1,2}}$ is an
$i$-disk with center $y$ and radius $\delta_{1,2}$, then
$D^i_{\delta_{1,2}}\subset B_{\delta_{1,1}}(x)$.

\; For $z\in B_{\delta_{1,3}}(x)$ and $l^{c+}_{\delta_{1,2}}(z)$
is a central curve at the right of $z$ with length $\delta_{1,2}$
and $z$ is one of its extreme points, suppose
$l^{c-}_{\delta_{1,2}}(z)$ is a central curve at the left of $z$
with length $\delta_{1,2}$ and $z$ is one of its extreme points,
let $l^c_{\delta_{1,2}}(z)=l^{c+}_{\delta_{1,2}}(z)\bigcup
l^{c-}_{\delta_{1,2}}(z)$, then $\#\{l^c_{\delta_{1,2}}(z)\bigcap
W^{uu}_{\delta_{1,2}}(y)\}=1$ and they are transverse intersect.
Suppose $z \notin W^{uu}_{\delta_{1,2}}(y)$, then if
$l^{c+}_{\delta_{1,2}}\bigcap W^{uu}_{\delta_{1,2}}(y)\neq \phi$,
we say $z$ is at $x$-left of $y$; if $l^{c-}_{\delta_{1,2}}\bigcap
W^{uu}_{\delta_{1,2}}(y)\neq \phi$, we say $z$ is at $x$-right of
$y$. It's easy to show when $z$ is at $x$-right of $y$, it's not
at $x$-right of $y$ anymore.

\; For simplicity, we just call $z$ at the left of
$W^{uu}_{loc}(y)$ or the right of $W^{uu}_{loc}(y)$.

\item[P6)] For any $x\in \overline{U_1}$, any
$\delta<\delta_{1,2}$, there exists $\delta^*\ll \delta$ such that
for $y \in B_{\delta^*}(x)\bigcap \Lambda_1$, if we have $z \in
B_{\delta^*}(x)\bigcap \Lambda_1$ also, then
$\#\{l^c_\delta(z)\bigcap W^{uu}_{\delta_{1,2}}(y)\}=1$ and they
are transverse intersect ($l^c_\delta(z)$ is defined in P5).

\item[P7)] For any $0<\delta^*<2\delta_{1,1}$, there exists a
$\delta^{**}$such that if $\Gamma$ is a central curve in
$\overline{U_1}$ with $length(\Gamma)<2\delta_{1,1}$, for $x,y\in
\Gamma$ and suppose the segment in $\Gamma$ connecting $x$ and $y$
has length bigger than $\delta^*$, then $d(x,y)>\delta^{**}$.

\item[P8)] For any $x\in \overline{U_1}$, any central curve $l$ in
$B_{\delta_{1,2}}(x)$ will have length smaller than
$\delta_{1,1}$.

\; For $y\in B_{\delta_{1,2}}(x)\bigcap \Lambda^+_1$, we can let
$W^{uu}_{\delta_{1,1}}(y)\bigcap B_{\delta_{1,2}}(x)$ always just
have one connected components, and $W^{uu}_{\delta_{1,1}/2}(y)$
divides $B_{\delta_{1,2}}(x)$ into two connected components: the
left part and the right part.

\; If $z_1,z_2\in B_{\delta_{1,2}}(x)$ are on the different side
of $B_{\delta_{1,2}}(x)\bigcap W^{uu}_{\delta_{1,1}/2}(y)$ and
there is a central curve $l\subset B_{\delta_{1,2}}(x)$ connecting
them, then $\#\{l\bigcap W^{uu}_{\delta_{1,1}/2}(y)\}=1$.

\item[P9)]Let $x\in \overline{U_1}$, suppose $y_1,y_2\in
B_{\delta_{1,2}}(x)\bigcap \Lambda_1^+$ and there exists a central
curve $l$ in $B_{\delta_{1,2}}(x)$ connects them, so by P8)
$length(l)<\delta_{1,1}$, now we know
$W^{uu}_{\delta_{1,1}/2}(y_1)\bigcap
W^{uu}_{\delta_{1,1}/2}(y_2)=\phi$ (otherwise $y_1\in
W^{uu}_{\delta_{1,1}}(y_2)$, then $\#\{l\bigcap
W^{uu}_{\delta_{1,1}}(y_1)\}\geq 2$, it contradicts with P4), it
means $W^{uu}_{\delta_{1,1}/2}(y_1)$ and
$W^{uu}_{\delta_{1,1}/2}(y_2)$ divide $B_{\delta_{1,2}}(x)$ into
three connected components. Suppose $y_1$ is at $x$-left of $y_2$,
then for any point $z\in \Lambda^+_1$ which are on the left of
$W^{uu}_{\delta_{1,1}/2}(y_2)\bigcap B_{\delta_{1,2}}(x)$ and on
the right of $W^{uu}_{\delta_{1,1}/2}(y_1)\bigcap
B_{\delta_{1,2}}(x)$, we have $W^{uu}_{\delta_{1,1}/2}(z)\bigcap
W^{uu}_{\delta_{1,1}/2}(y_i)=\phi\ _{(i=1,2)}$ and
$W^{uu}_{\delta_{11}/2}(z)\bigcap l \neq \phi$.

\item[P10)] A $C^1$ curve $\Gamma$ in $\overline{U_1}$ is called a
central segment if $f^i(\Gamma)\subset \overline{U_1}$ for all
$i\in \mathbb{Z}$ and it always tangents with $C^c_{a_0}$. Then
$\Gamma \subset \Lambda_1$ and it's easy to know that for any $x
\in \Gamma$, we have $T_x\Gamma =\widetilde{E^c_1}(x)$. On
$\Gamma$ we have normally hyperbolic splitting
$\widetilde{E^c_1}\oplus \widetilde{E^u}|_{\Gamma}$ since
$T_x\Gamma =\widetilde{E^c_1}(x)$, by the property of normally
hyperbolic manifold, $\bigcup \limits_{x\in
\Gamma}W^{uu}_{\delta_{1,1}/2}(x)$ is a submanifold ($dim=d$) with
boundary, we denote it $W^u_{\delta_{1,1}/2}(\Gamma)$.

\item[P11)] For any $\varepsilon>0$, if we have a family of
central segment $\{\Gamma_n\}^\infty_{n=1}$ with
$length(\Gamma_n)>\varepsilon$, there exists $\delta>0$ such that
$vol(W^u_{\delta_{1,1}/2}(\Gamma_n))>\delta$, so we can find
$n_i\neq n_j$ such that $W^u_{\delta_{1,1}/2}(\Gamma_{n_i})\bigcap
W^u_{\delta_{1,1}/2}(\Gamma_{n_j})$ $\neq \phi$.
\end{itemize}
\end{lem}

\subsection{Crovisier's central model}

In this subsection, let's fix
$U,U_1,U_2,\Lambda_1,\delta_0/2>\delta_{1,1}>\delta_{1,2}>\delta_{1,3}>0,$
and $a_0$ given by lemma \ref{5.3}, we'll introduce Crovisier's
central model. By his work, we can get some dynamical property for
the central curve $W^c_{\delta_{1,1}}(x)$ where $x\in \Lambda_1$.
The main result in this subsection is lemma $5.11$.

\begin{dfn}\label{5.4} A central model is a pair ($\widetilde{K},\widetilde{f}$) where
\begin{itemize}
\item[a)]  $\widetilde{K}$ is a compact metric space called the
base of the central model.

\item[b)] $\widetilde{f}$ is a continuous map from
$\widetilde{K}\times [0,1]$ into $\widetilde{K}\times [0,\infty)$

\item[c)]
$\widetilde{f}(\widetilde{K}\times\{0\})=\widetilde{K}\times\{0\}$

\item[d)] $\widetilde{f}$ is a local homeomorphism in a
neighborhood of $\widetilde{K}\times\{0\}$ : there exists a
continuous map $g:\ \widetilde{K}\times[0,1]\longrightarrow
\widetilde{K}\times [0,\infty)$ such that
$\widetilde{f}\circ\widetilde{g}$ and
$\widetilde{g}\circ\widetilde{f}$ are identity maps on
$\widetilde{g}^{-1}(\widetilde{K}\times[0,1])$ and
$\widetilde{f}^{-1}(\widetilde{K}\times[0,1])$ respectively.

\item[e)]  $\widetilde{f}$ is a skew product: there exits two map
$\widetilde{f}_1:\ \widetilde{K}\longrightarrow \widetilde{K}$ and
$\widetilde{f}_2:\ \widetilde{K}\times [0,1]\longrightarrow
[0,\infty)$ respectively such that for any $(x,t)\in
\widetilde{K}\times [0,1]$, one has
$\widetilde{f}(x,t)=(\widetilde{f}_1(x),\widetilde{f}_2(x,t))$.
\end{itemize}
$f$ general doesn't preserve $\widetilde{K}\times [0,1]$,
so the dynamics outside $\widetilde{K}\times \{0\}$ is only
partially defined.
\end{dfn}

The central model $(\widetilde{K},\widetilde{f})$ has a chain
recurrent central segment if it contains a segment $I=\{x\}\times
[0,a]$ contained in a chain recurrent class of
$f|_{\widetilde{K}\times [0,1]}$.

A subset $S\subset \widetilde{K}\times [0,1]$ of a product
$\widetilde{K}\times [0,\infty)$ is a strip if for any $x\in
\widetilde{K}$, the intersection $S\bigcap \{x\}\times [0,\infty)$
is a non-trivial interval.

In his remarkable paper \cite{1C2}, Crovisier got the following
important result.

\begin{lem}\label{5.5}(\cite{1C2} Proposition 2.5) Let $(\widetilde{K},\widetilde{f})$ be a central model with a
chain transitive base, then the two following properties are
equivalent:
\begin{itemize}
\item[a)]  There is no chain recurrent central segment.

\item[b)] There exists some strip $S$ in $\widetilde{K}\times
[0,1]$ that is arbitrarily small neighborhood of
$\widetilde{K}\times\{0\}$ and it's a trapping region for
$\widetilde{f}$ or $\widetilde{f}^{-1}$ : either
$\widetilde{f}(Cl(S))\subset Int(S)$ or
$\widetilde{f}^{-1}(Cl(S))\subset Int(S)$.
\end{itemize}
\end{lem}

\begin{rk} \label{5.6} If the central model $(\widetilde{K},\widetilde{f})$ has a chain recurrent
central segment and $\widetilde{K}\times \{0\}$ is transitive,
from Crovisier's proof, we can know for any small neighborhood $V$
of $\widetilde{K}\times \{0\}$, there exists a segment $x\times
[0,a]_{a\neq 0}$ contained in the same chain recurrent class of
$\widetilde{f}|_{V}$ with $\widetilde{K}\times\{0\}$.
\end{rk}

  An open strip $S\subset \widetilde{f}\times[0,1]$ satisfying $\widetilde{f}(Cl(S))\subset Int(S)$
or $\widetilde{f}^{-1}(Cl(S))\subset Int(S)$ will be called a
trapping strip.

\begin{dfn} \label{5.7} Let $f$ be a diffeomorphism of a manifold $M$,
$\Lambda,\Lambda_1,U,U_0,U_1,U_2,a_0,\delta_0/2>\delta_{1,1}>\delta_{1,2}>\delta_{1,3}>0$
are given in \S 5.1, where $\Lambda_1$ is a partial hyperbolic
invariant compact set of $f$ having a 1-dimensional central
bundle. A central model $(\widetilde{\Lambda}_1,\widetilde{f})$ is
a central model for $(\Lambda_1,f)$ if there exists a continuous
map $\pi:\ \widetilde{\Lambda}_1\times [0,\infty)\longrightarrow
M$ such that:
\begin{itemize}
\item[a)]  $\pi$ semi-conjugate $\widetilde{f}$ and $f$ :
$f\circ\pi=\pi\circ\widetilde{f}$ on
$\widetilde{\Lambda}_1\times[0,1]$

\item[b)] $\pi(\widetilde{\Lambda}_1\times\{0\})=\Lambda_1$

\item[c)]  The collection of map
$t\longrightarrow\pi(\widetilde{x},t)$ is a continuous family of
$C^1$ embedding of $[0,\infty)$ into $M$, parameterized by
$\widetilde{x}\in \widetilde{\Lambda_1}$.

\item[d)]  For any $\widetilde{x}\in \widetilde{\Lambda_1}$, the
curve $\pi(\widetilde{x},[0,\infty))\subset U$ has length bigger
than $\delta_{1,2}$ but smaller than $2\delta_{1,1}$, it's tangent
at the point $x=\pi(\widetilde{x},0)\in \Lambda_1$ to the central
bundle and it's a central curve ( that means the curve
$\pi(\widetilde{x},[0,\infty))$ tangents with the central cone
$C^c_{a_0}$).
\end{itemize}
\end{dfn}

\begin{rk} \label{5.8} From now, if $(\widetilde{\Lambda}_1,\widetilde{f})$ is a central model for ($\Lambda_1,f$)
and $\pi$ is the projection map, we'll denote the central model as
$(\widetilde{\Lambda}_1,\widetilde{f},\pi)$. Here I should notice
the reader that $\pi$ in this paper has two different meanings,
one denote the period of periodic point and another denote the
projection map of central model. If there is any confusion, I'll
point out.
\end{rk}

The following lemma shows the relation between central model and a
set with codimension-1 partial hyperbolic splitting.

\begin{lem} \label{5.9}([Cr2]) $\Lambda,\Lambda_1,U,U_1$ are given in $\S$5.1, then there exists a central
model $(\widetilde{\Lambda}_1,\widetilde{f},\pi)$ for
$(\Lambda_1,f)$. Let's denote $\widetilde{\Lambda}\subset
\widetilde{\Lambda}_1$ which satisfies $\pi^{-1}(\Lambda)\bigcap
(\widetilde{\Lambda}_1\times \{0\})=\widetilde{\Lambda}\times
\{0\}$, then $(\widetilde{\Lambda},\widetilde{f},\pi)$ is a
central model for $(\Lambda,f)$, and $\widetilde{\Lambda}\times
\{0\}$ is minimal.
\end{lem}

\begin{rk}\label{5.10}
\begin{itemize}
\item[1)] When the cental bundle $\widetilde{E^c_1}(\Lambda_1)$
has an $f$-orientation ( it means that
$\widetilde{E^c_1}(\Lambda_1)$ is orientable and $Df$ preserves
such orientation), we call the orientation 'right', then we can
get two central models
$(\widetilde{\Lambda^+_1},\widetilde{f}^+,\pi^+)$ and
$(\widetilde{\Lambda^-_1},\widetilde{f}^-,\pi^-)$ for
$(\Lambda_1,f)$, we call them the right model and the left model,
where $\pi^i\ _{(i=+,-)}$ is a bijection between
$\widetilde{\Lambda}^i_1\times \{0\}$ and $\Lambda_1$, and for
$\widetilde{x}^i\in \widetilde{\Lambda}^i_1$,
$\pi(\widetilde{x}^i\times[0,\infty))$ is a half of central curve
at the right ($i=+$) or left ($i=-$) of
$x=\pi(\widetilde{x}^i\times \{0\})$.

\item[2)] If $f$ doesn't preserve any orientation of
$\widetilde{E^c_1}(\Lambda_1)$, then $\pi:\
\widetilde{\Lambda}_1\longrightarrow \Lambda_1$ is two-one: any
point $x\in \Lambda_1$ has two preimages $\widetilde{x}^-$ and
$\widetilde{x}^+$ in $\widetilde{\Lambda}_1$, the homeomorphism
$\sigma$ of $\widetilde{\Lambda}_1$ which exchanges the preimages
$\widetilde{x}^+$ and $\widetilde{x}^-$ of any point $x\in
\Lambda_1$ commutes with $\widetilde{f}$.

In $\S$ 5.1, we know any point $x \in \Lambda_1$ has a local
orientation, then $\pi(\widetilde{x}^+\times [0,\infty))$ is a
central curve on the right of $x$, $\pi(\widetilde{x}^-\times
[0,\infty))$ is on the left of $x$, the union of them is a central
curve with central at $x$ and radius $\delta_{1,1}$.
\end{itemize}
\end{rk}

The following lemma is the main result in this subsection, it's
similar with [Cr]'s proposition 3.6, but a little stronger.

\begin{lem} \label{5.11} $f\in R$, $\Lambda$ is a non-trivial minimal set with a
codimension-1 partial hyperbolic splitting $E^c_1\oplus E^u$ where
$dim(E^c_1(\Lambda))=1$ and $E^c_1(\Lambda)$ is not hyperbolic.
Let $U,U_1,\Lambda_1$ be given in $\S$5.1, by lemma \ref{5.9},
$(\Lambda_1,f)$ has a central model
$(\widetilde{\Lambda}_1,\widetilde{f},\pi)$, then we can choose
$U_1$ properly such that
\begin{itemize}

\item[a)] either $(\widetilde{\Lambda}_1,\widetilde{f},\pi)$ has a
trapping region,

\item[b)] or $\Lambda$ is contained in a homoclinic class $C$, $C$
contains periodic points with index 1 and it's an index 0 fundamental
limit.
\end{itemize}
\end{lem}

\Pf: Let $\widetilde{\Lambda}\subset \widetilde{\Lambda}_1$
satisfy $\widetilde{\Lambda}\times \{0\}=\pi^{-1}(\Lambda)\bigcap
\widetilde{\Lambda}_1\times\{0\}$, then
$(\widetilde{\Lambda},\widetilde{f},\pi)$ is a central model for
$(\Lambda,f)$. Since now, we just denote
$\widetilde{\Lambda}\times\{0\}$ by $\widetilde{\Lambda}$.

At first, let's suppose $(\widetilde{\Lambda},\widetilde{f},\pi)$
has no trapping region, then by remark \ref{5.6}, for any small
neighborhood $V$ of $\widetilde{\Lambda}$ in
$\widetilde{\Lambda}\times [0,1]$, there exists a chain recurrent
central segment $x\times I$ in $V$ respecting the map
$\widetilde{f}$. By Crovisier's result ([Cr], proposition 3.6),
there exits a family of periodic points $\{p_n\}$ such that they
all belong to the same chain recurrent class with $\Lambda$ and
$\lim \limits_{n\rightarrow \infty} Orb(p_n)=\Lambda$, so $\Lambda
\subset H(p_n,f)_{n \geq 1}$. When $n$ is big enough,
$Orb(p_n)\subset \Lambda_1$, so $Orb(p_n)$ has a codimension-1
partial hyperbolic splitting $\widetilde{E}^c_1 \oplus
\widetilde{E}^u|_{Orb(p_n)}$, that means $p_n$ is an index 1
periodic point.

Now we claim that $H(p_n,f)$ is an index 0 fundamental limit.

\noindent{\bf Proof of the claim}: The argument is exactly the
same with the case $C$ in the proof of lemma \ref{4.3}, so here we
just give a sketch of the proof, we divide the proof to two cases.
\begin{itemize}

\item[A)]: there exists $\delta>0$ such that for any $p_n$, we
have $|Df^{\pi(p_n)}|_{\widetilde{E^c_1}(p_n)}|<e^{-\delta\pi(p_n)}$.

\item[B)]: for any $\frac{1}{m}>0$, there exists $p_{n_m}$ such
that $|Df^{\pi(p_{n_m})}|_{\widetilde{E^c_1}(p_{n_m})}|>e^{-\frac{1}{m}\pi(p_{n_m})}$.
\end{itemize}

In the first case, we use weakly selecting lemma, in case B, we
use lemma \ref{4.25}.          \qed

Now we suppose $(\widetilde{\Lambda},\widetilde{f},\pi)$ has a
trapping region $S$, we can suppose $\widetilde{f}(Cl(s))\subset
Int(S)$ always. Choose $\widetilde{\Lambda}_2$ an open
neighborhood of $\widetilde{\Lambda}$ in $\widetilde{\Lambda}_1$
small enough, we can get an open strip $S_2$ for
$\widetilde{\Lambda}_2$ (here open respect
$\widetilde{\Lambda}_2\times [0,1]$) such that:
\begin{itemize}
\item[a)] for any $\widetilde{x}\in \widetilde{\Lambda}$,
$\widetilde{x}\times [0,1]\bigcap S=\widetilde{x}\times
[0,1]\bigcap S_2$,

\item[b)] for any $\widetilde{x}\in \widetilde{\Lambda}_2$ and
$\widetilde{f}(\widetilde{x})\in \widetilde{\Lambda}_2$, we have
$\widetilde{f}(Cl((\widetilde{x}\times [0,1])\bigcap S_2))\subset
(\widetilde{f}(\widetilde{x})\times[0,1])\bigcap S_2$.
\end{itemize}

Choose $U^*$ neighborhood of $\Lambda$ small enough, let
$\Lambda^*=\bigcap \limits^{\infty}_{-\infty}f^i(\overline{U}^*)$,
then $\Lambda^* \subset \Lambda_1$, let
$\widetilde{\Lambda}^*\subset \widetilde{\Lambda}_1$ satisfies
$\widetilde{\Lambda}^*=\pi^{-1}(\Lambda^*)\bigcap
\widetilde{\Lambda}_1$, we'll have $\widetilde{\Lambda}^*\subset
\widetilde{\Lambda}_2$. Then consider the central model
$(\widetilde{\Lambda}^*,\widetilde{f},\pi)$ for $(\Lambda^*,f)$,
$S_2\bigcap (\widetilde{\Lambda}^*\times[0,1])$ is a trapping
region for $(\widetilde{\Lambda}^*,\widetilde{f},\pi)$.

  Now replace $U_1$ by $U^*$ and $\Lambda_1$ by $\Lambda^*$, we get a trapping region for
  $(\widetilde{\Lambda}_1,\widetilde{f},\pi)$.      \qed

\subsection{Proof of lemma 5.1}

Now we suppose $\Lambda$ is a non-trivial minimal set with a
codimension-1 partial hyperbolic splitting $E^c_1\oplus E^u$ where
$dim(E^c_1)=1$ and $E^c_1(\Lambda)$ is not hyperbolic. We divide
the proof of lemma 5.1 into two cases: $E^c_1(\Lambda)$ has an
$f$-orientation or not.\\

\noindent{\bf Proof of lemma 5.1} ( $E^c_1(\Lambda)$ has an
$f$-orientation)

Let $U_0$ be the small neighborhood of $\Lambda$ given in $\S$5.1
such that we can extend the splitting $E^c_1\oplus E^u|_\Lambda$
to $\overline{U}_0$, we denote the splitting
$T_xM=\widetilde{E}^c_1\oplus\widetilde{E}^u\ (x\in
\widetilde{U}_0)$. Suppose $U$ is any small neighborhood of
$\Lambda$ such that $\overline{U}\subset U_0$, then from lemma
\ref{5.3}, we can get open sets $U_2,U_1$ and $\Lambda_1=\bigcap
\limits_{i=-\infty}^{\infty}f^i(\overline{U}_1)$,
$a_0>0,0<\delta_{1,3}<\delta_{1,2}<\delta_{1,1}<\delta_0/2$ such
that they satisfy properties P1-P11 of lemma \ref{5.3} there.

  Since $E^c_1(\Lambda)$ has an $f$-orientation, $\widetilde{E}^c_1(\Lambda_1)$ has an $f$-orientation also, by
remark \ref{5.10} we get two central models: the right central
model $(\widetilde{\Lambda}^+_1,\widetilde{f}^+,\pi^+)$ and the
left central model
$(\widetilde{\Lambda}^-_1,\widetilde{f}^-,\pi^-)$, where for any
$\widetilde{x}^+\in \widetilde{\Lambda}^+_1$,
$\pi^+(\widetilde{x}^+\times [0,\infty))$ is a central curve at
the right of $x=\pi^+(\widetilde{x}^+\times \{0\})$ and
$\delta_{1,2}<length(\pi^+(\widetilde{x}^+\times
[0,\infty)))<2\delta_{1,1}$, so $\pi^+(\widetilde{x}^+\times
[0,\infty))\subset B_{2\delta_{1,1}}(x)\subset U$. For any
$\widetilde{x}^-\in \widetilde{\Lambda}^-$, we have the similar
property.

At first, we consider the right central model
$(\widetilde{\Lambda}^+_1,\widetilde{f}^+,\pi^+)$, if the right
central model doesn't have trapping region, by lemma \ref{5.11},
$\Lambda$ is contained in a homoclinic class $H(p,f)$ which
contains an index 1 periodic point and the homoclinic class is an
index 0 fundamental limit, then we've proved lemma \ref{5.1}, so
now we suppose that there exists a trapping region $S^+$ for the
right central model. By the similar argument for the left central
model, we can suppose it has a trapping region $S^-$ also.\\

\noindent{\bf Claim}: $\Lambda$ is an index 0 fundamental limit.\\

\noindent {\bf Proof of the claim:} If $\Lambda$ is not an index 0
fundamental limit, since $\Lambda$ has a codimension-1 dominated
splitting, $\Lambda$ should be an index 1 fundamental limit. By
generic property 5 of proposition \ref{3.1}, there exists a family
of index 1 periodic points $\{p_n\}$ such that $\lim
\limits_{n\rightarrow \infty}Orb(p_n)=\Lambda$ and they are index
stable, then by Gan's lemma, there exists a subsequence of
periodic points $\{p_{n_m}\}^\infty_{m=1}$ in $C$. Now with the
same argument of the case $C$ in the proof of lemma \ref{4.3}, we
can show $\Lambda$ satisfies weakly selecting lemma, by weakly
selecting lemma \ref{4.21}, $\Lambda$ is an index 0 fundamental
limit, that's a contradiction. \qed

Since $\Lambda$ is an index 0 fundamental limit, by generic
property 5) of proposition \ref{3.1}, there exists a family of
sources $\{p_n\}^\infty_{n=1}$ of $f$ satisfying $\lim
\limits_{n\rightarrow \infty}Orb(p_n)=\Lambda$. We can suppose
$Orb(p_n)\subset U_2$ always and let $\widetilde{p}^i_n\in
\widetilde{\Lambda}^i_1\ _{(i=+,-)}$ such that
$\pi^{(i)}(\widetilde{p}^i_n\times \{0\})=p_n$, then
$(\widetilde{f}^i)^{\pi(p_n)}(\widetilde{p}^i_n)=\widetilde{p}^i_n$.
Denote $\widetilde{p}_n^{+(-)}\times
I_n^{+(-)}=(\widetilde{p}^{+(-)}_n\times [0,\infty))\bigcap
S^{+(-)}$ and
$\gamma^{+(-)}_n=\pi^{+(-)}(\widetilde{p}^{+(-)}\times
I^{+(-)}_n)$, let $\gamma_=\gamma^+_n\bigcup \gamma_n^-$, then
$\gamma_n$ is a central curve with center at $p_n$. Since
$length(\gamma_n^{+(-)})<2\delta_{1,1}$, we have $\gamma_n\subset
B_{2\delta_{1,1}}(p_n)\subset U_1$.

We've suppose $S^{\pm}$ is a trapping region, then
$\widetilde{f}^{+(-)}(\overline{{S}^{+(-)}})\subset
Int({S}^{+(-)})$ or
$(\widetilde{f}^{+(-)})^{-1}(\overline{{S}^{+(-)}})\subset
Int({S}^{+(-)})$. In the first case, we say the trapping region is
1-step contracting, in the second case we say it's 1-step
expanding. When $S^i$ is 1-step contracting case, we have
$(\widetilde{f}^i)^{\pi(p_n)}(\widetilde{p}^i_n\times
\overline{I}^i_n)\subset \widetilde{p}^i_n\times I^i_n$, so
$f^{\pi(p_n)}(\overline{\gamma^i_n})\subset \gamma^i_n$ for
$i=+,-$ and there exists $\delta>0$ doesn't depend on $n$ such
that $length(\gamma^i_n\setminus
f^{\pi(p_n)}(\overline{\gamma^i_n}))>\delta$ for all $n\geq 1$. If
$S^i$ is 1-step expanding, we'll still have
$length(\gamma^i_n\setminus
f^{-\pi(p_n)}(\overline{\gamma^i_n}))>\delta$ for all $n\geq 1$.

Since $\gamma^i_n$ is either expanding or contracting for
$f^{\pi(p_n)}$, let $\Gamma^i_n=\bigcap \limits_{j=-\infty}^\infty
f^{j\pi(p_n)}(\gamma^i_n)\ _{(i=+,-)}$, we'll have
$f^{\pi(p_n)}(\Gamma^i_n)=\Gamma^i_n\ _{(i=+,-)}$ where
$\Gamma^i_n$'s extreme points are periodic points. When
$\Gamma^i_n$ is not trivial, we denote $q_n^i\ _{(i=+,-)}$ the
extreme periodic point different with $p_n$, if $\Gamma^i_n$ is
trivial, we just let $q_n^i=p_n$. We let
$\Gamma_n=\Gamma^+_n\bigcup\Gamma^-_n$ and
$h^i_n=\gamma_n^i\setminus\Gamma^n_i\ _{(i=+,-)}$, then
$\Gamma_n\subset \Lambda_1$, $h^i_n\subset U_1$. It's easy to know
that $h^i_n$ is in the stable (unstable) manifold of $q_n^i$ if
$S^i$ is 1-step contracting (expanding). And since $f$ is a
Kupka-Smale diffwomorphism, $f^{\pi(p_n)}|_{\Gamma_n}$ is also a
Kupka-Smale diffeomorphism and just has finite sinks and sources
(respect $f^{\pi(p_n)}|_{\Gamma_n}$).

\begin{lem}\label{5.12} If $\Gamma_n\bigcap\Gamma_m\neq \phi$,
 then $\Gamma_n\bigcap\Gamma_m$ is a connected central curve, and
$\Gamma_n\bigcup\Gamma_m$ is a central segment.
\end{lem}

\Pf: We need prove some lemmas at first.

\begin{lem} \label{5.13} let $x\in \Gamma_n\bigcap\Gamma_m$ and $x$
is not a periodic point, $x_1\in \Gamma_n$ is the nearest periodic
point at the left of $x$ and $x_2\in \Gamma_n$ is the nearest
periodic point at the right of $x$. Denote $I_n\subset \Gamma_n$
the segment connecting $x_1$ and $x_2$, then $I_n \subset
\Gamma_m$.
\end{lem}

\Pf: By the assumption, $f^{\pi(p_n)}$ has no any other fixed
point in $I_n$, so for $x_1$ and $x_2$, one of them is sink for
$f^{\pi(p_n)}|_{\Gamma_n}$ and another is source for
$f^{\pi(p_n)}|_{\Gamma_n}$. We suppose $x_1$ is the source, then
$\lim \limits_{i\rightarrow\infty}f^{i\pi(p_n)}(x)\longrightarrow
x_2$ and $\lim
\limits_{i\rightarrow\infty}f^{-i\pi(p_n)}(x)\longrightarrow x_1$.
Since $\Gamma_m$ is a periodic central segment with period
$\pi(p_m)$ and $x\in \Gamma_m$, we have
$f^{i\pi(p_n)\pi(p_m)}(x)\in \Gamma_m$ for all $i \in \mathbb{Z}$,
so $x_2=\lim
\limits_{i\rightarrow\infty}f^{i\pi(p_n)\pi(p_m)}(x)\in \Gamma_m$
and $x_1=\lim
\limits_{i\rightarrow\infty}f^{-i\pi(p_n)\pi(p_m)}(x)\in
\Gamma_m$.

  Now denote $I_m$ the central segment in $\Gamma_m$ connecting $x_1$ and $x_2$.

  We claim that $I_n=I_m$.\\

\noindent{\bf Proof of the claim}: If it's not true, there exists
$y\in Int(I_n)$, $z\in W^{uu}_{\delta_{1,1}}(y)\bigcap I_m$ and
$z\neq y$.

    For any $\varepsilon>0$, consider $a=f^{i\pi(p_n)\pi(p_m)}(y)$
where $i$ is very big, then $a\in I_n$ and it's near $x_2$ very
much. Let $b\in W^{uu}_{\delta_{1,1}}(a)\bigcap I_m$, recall that
$I_n$ and $I_m$ are tangent at $\widetilde{E^c_1}(x_2)$, when $i$
is big enough, there exists a curve $l$ in
$W^{uu}_{\delta_{1,1}}(a)$ connecting $a$ and $b$ with
$length(l)<\varepsilon$.

\begin{figure}[h]
\begin{center}
\psfrag{12}{$W^{uu}_{\delta_{1,1}}(x_1)$}\psfrag{13}{$W^{uu}_{\delta_{1,1}}(y)$}
\psfrag{14}{$W^{uu}_{\delta_{1,1}}(a)$}\psfrag{15}{$W^{uu}_{\delta_{1,1}}(x_2)$}
\psfrag{16}{$x_1$}\psfrag{17}{$y$}\psfrag{18}{$a$}\psfrag{19}{$x_2$}
\psfrag{20}{$z$}\psfrag{21}{$b$}\psfrag{22}{$I_n$}\psfrag{23}{$I_m$}\psfrag{x}{$x$}
\includegraphics[height=2.0in]{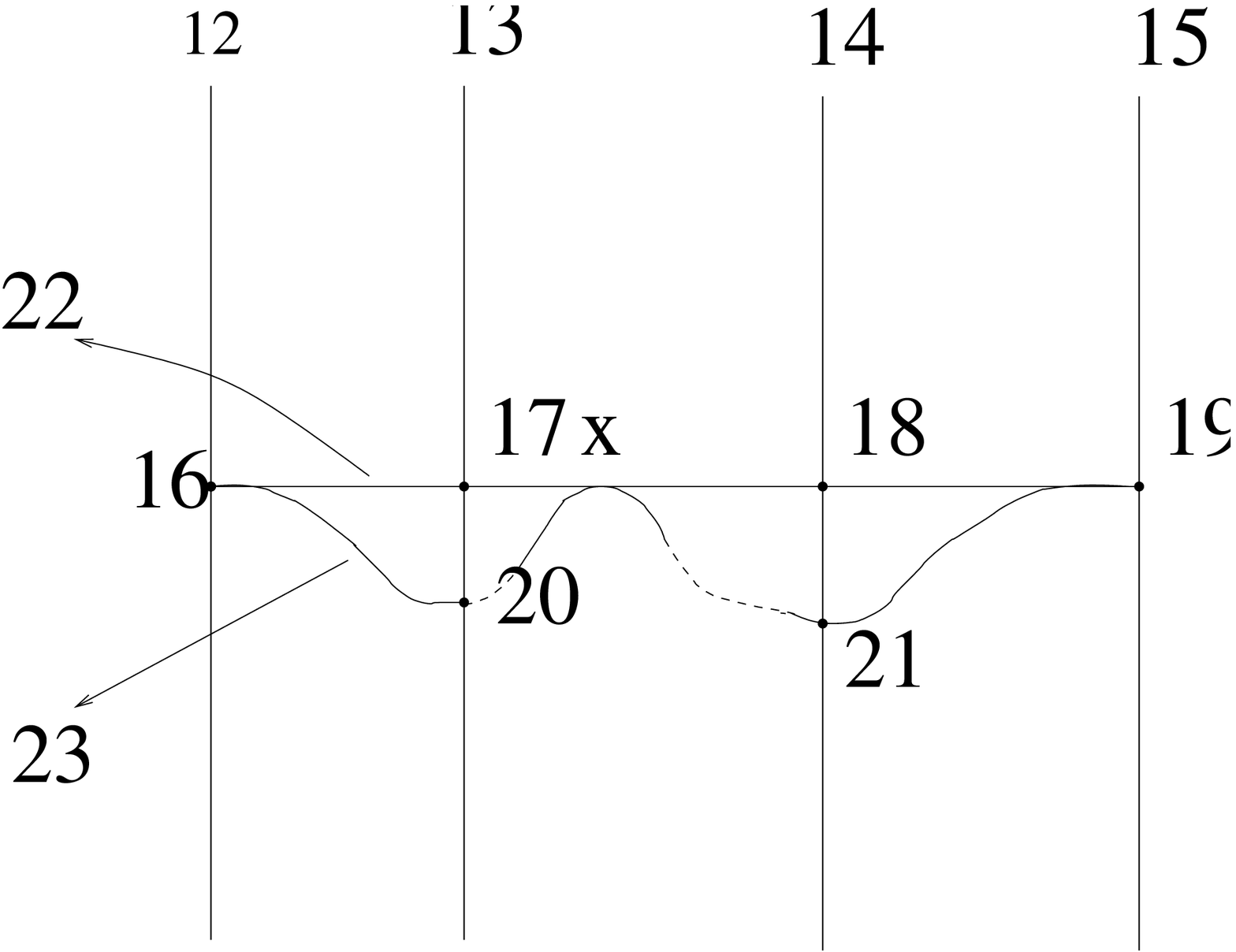}
\end{center}
\end{figure}

 Now it's easy to know $f^{-i\pi(p_n)\pi(p_m)}(b)\in W^{uu}_{\delta_{1,1}}(y)\bigcap \Gamma_m$.
  By P4 of lemma \ref{5.3}, $\#\{W^{uu}_{\delta_{1,1}}(y)\bigcap \Gamma_m\}=1$,
 so $f^{-i\pi(p_n)\pi(p_m)}(b)=z$, then $f^{-i\pi(p_n)\pi(p_m)}(l)$ is a curve connecting
 $y$ and $z$, by P3 of lemma \ref{5.3}, we'll have
 $length(f^{-i\pi(p_n)\pi(p_m)}(l))<\varepsilon\cdot \lambda^{i\pi(p_n)\pi(p_m)}$.

 Since $\varepsilon$ can be chosen
arbitrarily small, we get $y=z$, that's a contradiction. \qed

      By the claim, we finish the proof of lemma \ref{5.13}.
      \qed

    We still need the following result.

\begin{lem}\label{5.14} Let $x\in \Gamma_n\bigcap \Gamma_m$ and $x$
be a fixed point of $f^{\pi(p_n)}|_{\Gamma_n}$ and
$f^{\pi(p_m)}|_{\Gamma_m}$, suppose $\Gamma_n$ and $\Gamma_m$ both
have points on the right of $x$. Let $x_n\in \Gamma_n$ be the
nearest fixed point of $f^{\pi(p_n)}|_{\Gamma_n}$ on the right of
$x$ and $x_m\in \Gamma_m$ be the nearest fixed point of
$f^{\pi(p_m)}|_{\Gamma_m}$ on the right of $x$. Denote $I_n\subset
\Gamma_n$ the central segment in $\Gamma_n$ connecting $x$ and
$x_n$, $I_m\subset \Gamma_m$ the central segment in $\Gamma_m$
connecting $x$ and $x_m$, then $I_n=I_m$.
\end{lem}

\Pf:  At first, we claim that either
$W^{uu}_{\delta_{1,1}}(x_n)\bigcap I_m\neq \phi$ or
$W^{uu}_{\delta_{1,1}}(x_m)\bigcap I_n\neq \phi$.

\noindent{\bf Proof of the claim}: Suppose
$W^{uu}_{\delta_{1,1}}(x_n)\bigcap I_m\neq \phi$, we know that
$x_m$ is on the left of $W^{uu}_{\delta_{1,1}}(x_n)$, recall that
$x_m$ is on the right of $x$, so by P9 of lemma \ref{5.3},
$W^{uu}_{\delta_{1,1}}(x_m)\bigcap I_n\neq \phi$. \qed

Now we suppose $W^{uu}_{\delta_{1,1}}(x_n)\bigcap I_m=y\neq \phi$,
then $y\in I_m\setminus \{x\}$, it's easy to know
$f^{-i\pi(p_n)\pi(p_m)}(y)\in W^{uu}_{\delta_{1,1}}(x_n)\bigcap
I_m$ for $i\geq 1$, so $f^{-i\pi(p_n)\pi(p_m)}(y)=y$. But $\lim
\limits_{i\rightarrow \infty}
f^{-i\pi(p_n)\pi(p_m)}(y)\longrightarrow x_n$, so $x_n=y$. It
means that $x_n\in I_m\setminus \{x\}$, so $x_n=x_m$. By the same
argument in lemma \ref{5.13}, we can prove $I_n=I_m$. \qed

Now let's continue the proof of lemma \ref{5.12}.

Let $\Gamma=\Gamma_n\bigcap \Gamma_m$, $x\in \Gamma$ be the left
extreme point of $\Gamma$, then by lemma \ref{5.13}, $x$ should be
a periodic point and on the left of $x$, there doesn't contain
points of at least one of the segment $\Gamma_n$ or $\Gamma_m$.
Let $y\in \Gamma$ be the right extreme point of $\Gamma$, then on
the right of $y$, there doesn't contain points of at least one of
the segments $\Gamma_n$ or $\Gamma_m$.

When $x=y$, $\Gamma_n$ and $\Gamma_m$ are on different side of
$x$, $\Gamma_n\bigcup \Gamma_m$ is obviously a central segment.

When $x\neq y$, let $I$ be the maximal central curve in $\Gamma$
containing $x$, let $z$ be the right extreme point in $I$, by
lemma \ref{5.13}, $z$ should be a periodic point. If $z\neq y$,
$y$ is on the right of $z$ and $y\in \Gamma_n\bigcap \Gamma_m$, so
by lemma \ref{5.14}, $I$ will contain a central segment on the
right of $z$, that's a contradiction with the maximalicity of $I$,
so $z=y$. It means that $I=\Gamma_n\bigcap \Gamma_m$ is an
interval, and $x,y$ are its extreme points on the left and right,
and $\Gamma_n$ and $\Gamma_m$ can not both have points on the left
of $x$, they can not both have points on the right of $y$ also,
it's easy to see now that $\Gamma_n\bigcup \Gamma_m$ is a central
curve. \qed

Now we divide the proof of lemma \ref{5.1} to three cases
depending on the contracting or expanding properties of the two
central models.\\

\noindent{\bf Case A}: {\it Two central models have 1-step
expanding
properties.}\\

 In this case, for any $\gamma_n$, we have $f^{-i}(\gamma_n)\in U_1$ for $i \geq 1$, it means
$\gamma\subset \Lambda^+_1$, and any $x\in \gamma_n$ will have
uniform size of unstable manifold $W^{uu}_{\delta_{1,1}}(x)$. Let
$W^u_{\delta_{1,1}/2}(\gamma_n)=\bigcup \limits_{x\in
\gamma_n}W^{uu}_{\delta_{1,1}/2}(x)$, by the property of normally
hyperbolic submanifold, $W^u_{\delta_{1,1}/2}(\gamma_n)$ is a
submanifold ($dim=d$) with boundary, it's easy to know that
$W^u_{\delta_{1,1}/2}(\gamma_n)$ has uniform size, that means
there exists an $\varepsilon>0$ such that
$B_\varepsilon(p_n)\subset W^u_{\delta_{1,1}/2}(\gamma_n)$ for all
$n\geq 1$. Suppose $\lim \limits_{n\rightarrow \infty}p_n=p\in
\Lambda$, then when $n$ is big enough, $p\in
B_\varepsilon(p_n)\subset W^u_{\delta_{1,1}/2}(\gamma_n)$, so
$\lim \limits_{i\rightarrow
\infty}f^{-i\pi(p_n)}(p_n)\longrightarrow$ some periodic point
$z\in \Gamma_n$, it means $z\in \Lambda$. But $\Lambda$ is a
non-trivial minimal set of $f$, that's a contradiction. \qed

\noindent{\bf Case B}: {\it Left central model is 1-step
contracting and
the right central model is 1-step expanding.}\\

 Let's consider $\gamma^+_n$,
with the same argument in case A, it has uniform size of unstable
manifold $W^u_{\delta_{1,1}/2}(\gamma^+_n)=\bigcup \limits_{x\in
\gamma^+_n}W^{uu}_{\delta_{1,1}/2}(x)$ ( it's because
$length(\gamma^+_n)>length(h^+_n)>\delta$), so there exists an
$\varepsilon>0$ such that
$Vol(W^u_{\delta_{1,1}/2}(\gamma^+_n))>\varepsilon$.

 Now we claim that for any sequence $\{n_i\}^\infty_{i=1}$,
there exists $i_0$ and a sequence $i_0<i_1<i_2<\cdots$ such that
for any $j>0$, $W^u_{\delta_{1,1}/2}(\gamma^+_{n_{i_j}})\bigcap
W^u_{\delta_{1,1}/2}(\gamma^+_{n_{i_0}})\neq \phi$.\\

\noindent{\bf Proof of the claim}: Suppose that the claim is not
true, then we can find a subsequence $\{n_{i_j}\}^\infty_{j=1}$
such that $W^u_{\delta_{1,1}/2}(\gamma^+_{n_{i_{j_0}}})\bigcap
W^u_{\delta_{1,1}/2}(\gamma^+_{n_{i_j}})=\phi$ for $j_0\in
\mathbb{N}$ and $j>j_0$, it's a contradiction with
$Vol(W^u_{\delta_{1,1}/2}(\gamma^+_{n_i}))>\varepsilon$, since
we'll have $Vol(M)>\sum
\limits_jVol(W^u_{\delta_{1,1}/2}(\gamma^+_{n_{i_j}}))=\infty$.
\qed

By the above claim, we can find a subsequence
$\{n_i\}^\infty_{i=1}$ such that for any $i_0\in \mathbb{N}^+$, we
can get $W^u_{\delta_{1,1}/2}(\gamma^+_{n_i})\bigcap
W^u_{\delta_{1,1}/2}(\gamma^+_{n_{i_0}})\neq \phi$ for $i\geq
i_0$. Since $f$ is a Kupka-Smale diffeomorphism, on $\Gamma_{n_i}$
it just has finite periodic points. So when we fix $i_0$, we can
let $i$ big enough such that $p_{n_i}\notin \gamma_{n_{i_0}}$. It
means that we can choose a subsequence
$\{(\Gamma_{n_i},\Gamma_{m_i})\}^\infty_{i=0}$ such that
$p_{m_i}\notin \Gamma_{n_i}$,
$W^u_{\delta_{1,1}/2}(\gamma^+_{n_i})\bigcap
W^u_{\delta_{1,1}/2}(\gamma^+_{m_i})\neq \phi$ and $\lim
\limits_{i\rightarrow \infty}(p_{n_i})=\lim \limits_{i\rightarrow
\infty}(p_{m_i})=x_0$ for some $x_0\in \Lambda$.

Since $W^u_{\delta_{1,1}/2}(\gamma^+_{n_i})\bigcap
W^u_{\delta_{1,1}/2}(\gamma^+_{m_i})\neq \phi$, suppose $y_i\in
W^u_{\delta_{1,1}/2}(\gamma^+_{n_i})\bigcap
W^u_{\delta_{1,1}/2}(\gamma^+_{m_i})$, then $$\lim
\limits_{j\rightarrow
\infty}f^{-j\pi(p_{n_i})\pi(p_{m_i})}(y_i)\longrightarrow
\Gamma^+_{n_i}\ \text{and}\ \lim \limits_{j\rightarrow
\infty}f^{-j\pi(p_{n_i})\pi(p_{m_i})}(y_i)\longrightarrow
\Gamma^+_{m_i},$$ so $\Gamma^+_{n_i}\bigcap \Gamma^+_{m_i}\neq
\phi,$ by lemma \ref{5.12}, $\Gamma_{n_i}\bigcup \Gamma_{m_i}$ is
a central segment.

 For simplicity, we suppose $p_{m_i}$ is on the right
of $p_{n_i}$ for all $i\in \mathbb{N}$, the proof of the other
case is similar. Since $p_{m_i}\notin \Gamma_{n_i}$ and
$\Gamma_i=\Gamma_{n_i}\bigcup \Gamma_{m_i}$ is a central curve.
$p_{m_i}$ is on the right of $q^+_{n_i}$ also. Recall that
$q^+_{n_i}$ is a source for $f^{\pi(p_{n_i})}|_{\Gamma_{n_i}}$,
and $h^+_{n_i}$ belongs to its basin, so $h^+_{n_i}\bigcap
W^{uu}_{\delta_{1,1}/2}(p_{m_i})=\phi$.

\begin{rk}: We don't know $h^+_{n_i}\subset \Gamma_{m_i}$ here.\end{rk}

We know that $h^+_{n_i}$ is a central curve on the right of
$q^+_{n_i}$ with length bigger than $\delta$, by property P6 of
lemma \ref{5.3}, there exists a $\delta^*$ such that
$d(q^+_{n_i},p_{m_i})>\delta^*$.( Since if
$d(q^+_{n_i},p_{m_i})<\delta^*$, we have
$l^+_\delta(q^+_{n_i})\bigcap W^{uu}_{\delta_{1,1}/2}(p_n)\neq
\phi$ where $l^+_\delta(q^+_{n_i})$ is any central curve at the
right of $q^+_{n_i}$ with length $\delta$ and $q^+_{n_i}$ is the
left extreme point of it, with the fact that $p_{m_i}$ is on the
right of $q^+_{n_i}$, we'll have $h^+_{n}\bigcap
W^{uu}_{\delta_{1,1}/2}(p_{m_i})\neq \phi$, that's a contradiction
because $h^+_{n_i}\subset W^u(q^+_{n_i})$). So especially, in the
central segment $\Gamma_i$, the distance between $p_{n_i}$ and
$p_{m_i}$ is bigger than $\delta^*$. By property P7 of lemma
\ref{5.3}, there exists $\delta^{**}>0$ such that
$d(p_{n_i},p_{m_i})>\delta^{**}$, it's a contradiction with $\lim
\limits_{i\rightarrow \infty}(p_{n_i})=\lim \limits_{i\rightarrow
\infty}(p_{m_i})=x_0\in \Lambda$. \qed

\noindent{\bf Case C}: {\it The two central models have 1-step
contracting properties.}\\

  In this case, replace by a subsequence, we can suppose for
$\{\Gamma_n\}^\infty_{n=1}$, we have $p_n\notin \bigcup
\limits_{i<n}\Gamma_i$.

\begin{lem}\label{5.15} There exists $n_0$ big enough such that for any
$n_1,n_2>n_0$, $n_1\neq n_2$, we always have
$W^u_{\delta_{1,1}/2}(\Gamma_{n_1})\bigcap
W^u_{\delta_{1,1}/2}(\Gamma_{n_2})=\phi$.
\end{lem}

\Pf£ºSuppose the lemma is not true, then we can choose $n_1$ and
$n_2$ arbitrarily big and satisfying
$W^u_{\delta_{1,1}/2}(\Gamma_{n_1})\bigcap
W^u_{\delta_{1,1}/2}(\Gamma_{n_2})\neq \phi$, then it's easy to
know $\Gamma_{n_1}\bigcap \Gamma_{n_2}\neq \phi$ and
$\Gamma_{n_1}\bigcup \Gamma_{n_2}$ is a central curve. We can
suppose $n_2>n_1$, then by the assumption of
$\{\Gamma_n\}^\infty_{n=1}$, we have $p_{n_2}\notin \Gamma_{n_1}$.

  We just suppose $p_{n_2}$ is on the right of $p_{n_1}$, since
$\Gamma=\Gamma_{n_1}\bigcup \Gamma_{n_2}$ is a central curve and
$p_{n_2}\notin \Gamma_{n_1}$, we can know $p_{n_2}$ is on the
right of $q^+_{n_1}$ also, and $q^+_{n_1}\in \Gamma_{n_2}$.

  We know that there exists a $\delta>0$ such that $length(h^{+(-)}_n)>\delta$
for all $n\geq 1$. And for such $\delta$, by proposition P6 of
lemma \ref{5.3}, there exists $0<\delta^*\ll \delta$ such that for
any $x,y\in \Lambda_1$, if $d(x,y)<\delta^*$, we have
$\#\{W^{uu}_{\delta_{1,1}/2}(x)\bigcap l^c_\delta(y)\}=1$ where
$l^c_\delta(y)$ is a central curve with center $y$ and on the two
sides of $y$ both have length $\delta$.

  Suppose $x\in \Gamma_m$ is the nearest periodic point on the right
side of $q^+_{n_1}$, and let $I\subset \Gamma_m$ the central
segment in $\Gamma_m$ connecting $q^+_{n_1}$ and $x$.

  Now we claim that $length(I)>\delta^*$.\\

\noindent{\bf Proof of the claim}: If $length(I)\leq \delta^*$,
then $d(q^+_{n_1},x)\leq \delta^*$ also. By the facts that $x$ is
on the right of $q_{n_i}^+$ and $h^+_{n_1}$ is a central curve
with length bigger than $\delta$, we have $h^+_{n_1}\bigcap
W^{uu}_{\delta_{1,1}/2}(x)\neq \phi$. Then for any $y\in Int(I)$,
$W^{uu}_{\delta_{1,1}/2}(y)\bigcap h^+_{n_1}\neq \phi$.

It's easy to know $I\nsubseteq h^+_{n_1}$ since $h^+_{n_1}$
contains no periodic point, so there exists $z\in h^+_{n_1}$ such
that $W^{uu}_{\delta_{1,1}/2}(z)\bigcap Int(I)=y\neq z$.

\begin{figure}[h]
\begin{center}
\psfrag{1}[][]{$W^{uu}_{\delta_{1,1}/2}(q_{n_1}^+)$}
\psfrag{3}[][]{$W^{uu}_{\delta_{1,1}/2}(z)$}\psfrag{4}{$W^{uu}_{\delta_{1,1}/2}(x)$}
\psfrag{5}{$q_{n_1}^+$}\psfrag{6}{$z_i$}\psfrag{7}{$z$}\psfrag{8}{$h_{n_1}^+$}\psfrag{9}{$a_i$}\psfrag{10}{$y$}
\psfrag{11}{$x$}
\includegraphics[height=1.5in]{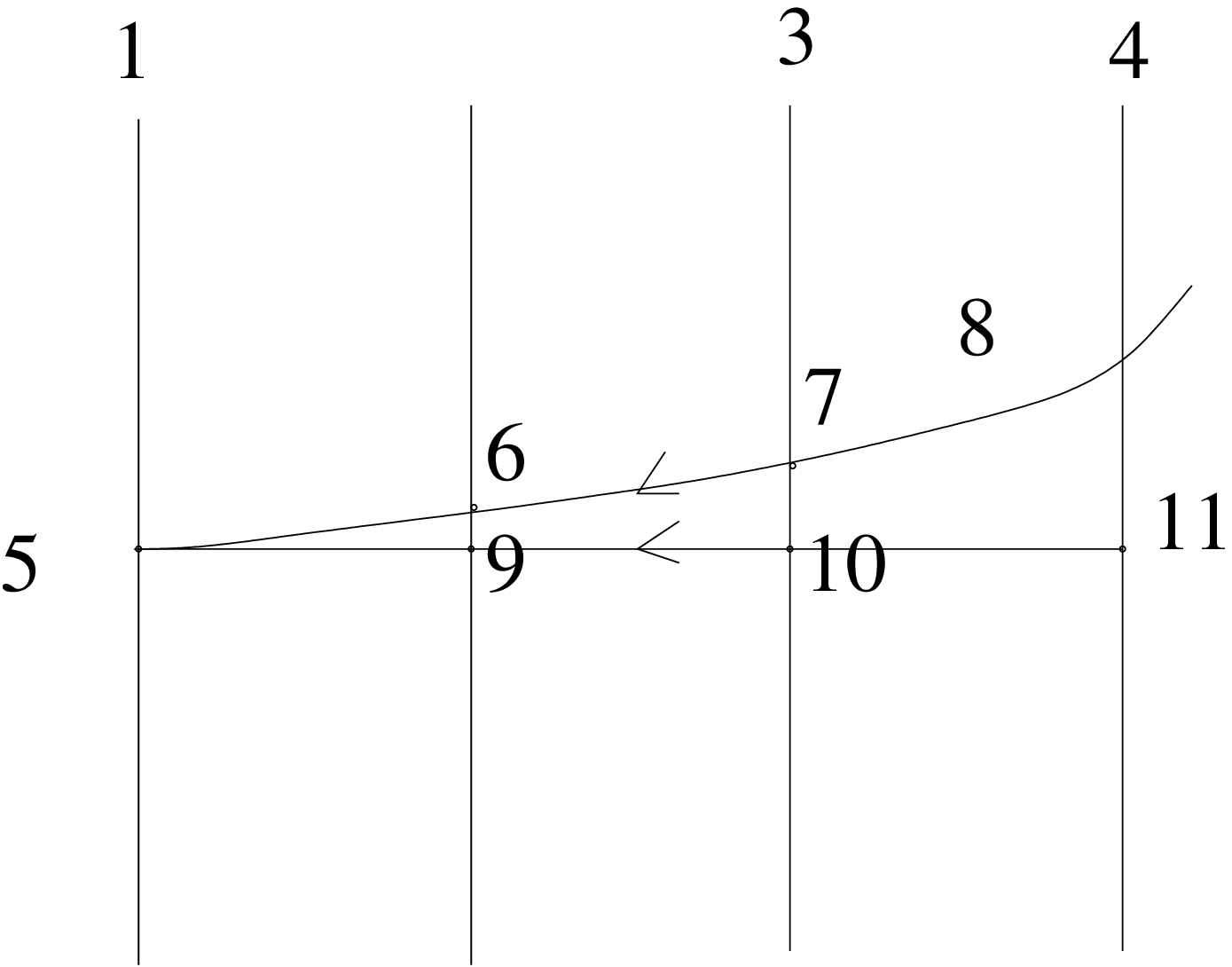}
\end{center}
\end{figure}

Because the two central models are 1-step contracting, $q^+_{n_1}$
is a sink for $f^{\pi(p_{n_1})}|_{\Gamma_{n_1}}$, then it's also a
sink for $f^{\pi(p_{n_1})\pi(p_{n_2})}|_\Gamma$ where
$\Gamma=\Gamma_{n_1}\bigcup \Gamma_{n_2}$. We can choose $i$ big
enough, such that $z_i=f^{i\pi(p_{n_1})\pi(p_{n_2})}(z)$ near
$q^+_{n_1}$ very much, let
$a_i=W^{uu}_{\delta_{1,1}/2}(z_i)\bigcap I$. Since $h^+_{n_1}$ and
$I$ are tangent at $q^+_{n_1}$ on $\widetilde{E^c_1}(q^+_{n_1})$,
for any $\varepsilon>0$, when $i$ big enough, there exists a curve
$l\subset W^{uu}_{\delta_{1,1}/2}(z_i)$ connecting $a_i$ and $z_i$
and $length(l)<\varepsilon$. Since
$f^{-i\pi(p_{n_1})\pi(p_{n_2})}(a_i)\in
W^{uu}_{\delta_{1,1}/2}(z)\bigcap I$, that means
$f^{-i\pi(p_{n_1})\pi(p_{n_2})}(a_i)=y$ and
$f^{-i\pi(p_{n_1})\pi(p_{n_2})}(l)$ is a curve connecting $z$ and
$y$. By property P3 of lemma \ref{5.3},
$length(f^{-i\pi(p_{n_1})\pi(p_{n_2})}(l))<\varepsilon\lambda^i$.
Since $i$ can be chosen arbitrarily big, we can get $y=z$, that's
a contradiction. \qed

Since $length(I)>\delta^*$, the segment in $\Gamma$ connecting
$p_{n_1}$ and $p_{n_2}$ will have length bigger than $\delta^*$
also, by property P7 of lemma \ref{5.3}, there exists
$\delta^{**}>0$ such that $d(p_{n_1},p_{n_2})>\delta^{**}$. But
recall that $\lim \limits_{n\rightarrow \infty}p_n\longrightarrow
x_0\in \Lambda$ and $n_1,n_2$ can be chosen arbitrarily big, we
can get $d(p_{n_1},p_{n_2})<\delta^{**}$, that's a contradiction.
\qed

With lemma \ref{5.15}, we can chosen $\{\Gamma_n\}^\infty_{n=1}$
such that if $n\neq m$, $W^u_{\delta_{1,1}/2}(\Gamma_n)\bigcap
W^u_{\delta_{1,1}/2}(\Gamma_n)=\phi$. Then by property P11 of
lemma \ref{5.3}, $\lim \limits_{n\rightarrow
\infty}length(\Gamma_n)=0$.

Choose $n_0$ big enough such that for $m \geq n_0$,
$d(p_m,p_{n_0})<\delta^*/4$ and $length(\Gamma_m)<\delta^*/4$, we
can suppose $p_m$ is on the right of $p_{n_0}$, then by
$W^u_{\delta_{1,1}/2}(\Gamma_n)\bigcap
W^u_{\delta_{1,1}/2}(\Gamma_n)=\phi$, we know that $p_m$ is on the
right of $q^+_{n_0}$ and $q^-_m$ is on the right of $q_{n_0}^+$
also.

Since $d(q^+_{n_0},q^-_m)\leq
d(q^+_{n_0},p_{n_0})+d(q^-_m,p_m)+d(p_{n_0},p_m)<length(\Gamma_{n_0})+\delta^*/4
+length(\Gamma_m)<\delta^*$, by Property P6 of lemma \ref{5.3} and
$length(h^+_{n_0})>\delta$, $length(h_m^-)>\delta$, we can get
$h^+_{n_0}\pitchfork W^{uu}_{\delta_{1,1}/2}(q^-_m)\neq \phi$ and
$h^+_m\pitchfork W^{uu}_{\delta_{1,1}/2}(q^+_{n_0})\neq \phi$.
Recall that $h^+_{n_0}\subset W^s(q^+_{n_0})$ and $h^-_m \subset
W^s(q^-_m)$, we can know $q^+_{n_0}$ and $q^-_m$ are in the same
homoclinic class.

When $m\longrightarrow \infty$, by
$length(\Gamma_m)\longrightarrow 0$ and $\lim
\limits_{m\rightarrow \infty}p_m\longrightarrow x_0\in \Lambda$,
we have $q^-_m\longrightarrow x_0$ also, so $x\in H(q^+_{n_0},f)$
and then $\Lambda \subset H(q^+_{n_0},f)$.

  Now we'll prove $H(q^+_{n_0},f)$ is an index 0 fundamental limit.

Recall that $Orb(q^+_{n_0})\subset U$ and $U$ can be chosen
arbitrarily small, so in fact we've proved that there exists a
family of periodic points $q_n$ with index 1 such that $\lim
\limits_{n\rightarrow \infty}Orb(q_n)=\Lambda$ and $\Lambda\subset
H(q_1,f)=H(q_2,f)=\cdots$.

By the same argument with case C in the proof of lemma \ref{4.3},
we can prove $H(q_1,f)$ is an index 0 fundamental limit. \qed

  Now let's keep on proving the other case of lemma \ref{5.1}.\\

\noindent{\bf Proof of lemma 5.1}($E^c_1(\Lambda)$ has no any
$f$-orientation):\\

In this case, we just have one central model, but locally we still
have orientation for $\widetilde{E^c_1}(\Lambda_1)$, and the two
sides have the same dynamical property: they are both 1-step
expanding or they are both 1-step contracting. All the other
argument is the same with the case where $E^c_1(\Lambda)$ has an
$f$-orientation. \qed

 \bt{99}
\bib{1ABC}
F. Abdenur, C. Bonatti, and S. Crovisier, Global dominated
splittings and the $C^1$ Newhouse phenomenon, {\it Proceedings of
the American Mathematical Society} {\bf 134}, (2006), 2229-2237.\
\bib{1ABCD}
F. Abdenur, C. Bonatti, S. Crovisier, L.J. Diaz and L. Wen, {\it
Periodic points and homoclinic classes,} {\it preprint} (2006).\
\bib{1AS}
R. Abraham and S. Smale, Nongenericity of $\Omega$
-stability, Global analysis I, {\it Proc. Symp. Pure
Math. AMS} {\bf 14} (1970), 5-8.\
\bib{1A }
M-C. Arnaud, Creation de connexions en topologie $C^1$. {\it
Ergodic Theory and Dynamical System}{\bf 31} (2001), 339-381.\
\bib{1BC}
C. Bonatti and S. Crovisier, Recurrence et genericite(French),
{\it Invent. math.,} {\bf 158} (2004), 33-104\
\bib{1BDP}
C. Bonatti L. Diaz and E. Pujals, A $C^1$-generic dichotomy for
diffeomorphisms: weak forms of hyperbolicity or infinitely many
class or sources, {\it Ann. of Math.} {\bf 158} (2003), 355-418\
\bib{1BD}
C. Bonatti, L. J. D\'{i}az, Connexions h\'{e}t\'{e}rocliniques et
g\'{e}n\'{e}ricite d'une infinit\'{e} de puits ou de sources, {\it
Annales Scientifiques de l'¨¦cole Normal Sup¨¦rieure de Paris},
{\bf 32 (4)}, (1999) 135-150,\
\bib{1BDV}
C. Bonatti, L. J. D\'{i}az, and M. Viana, Dyanamics beyond uniform
hyperbolic, {\bf Volume 102} of {\it Encyclopaedia of Mathematical
Sciences.} Springer- Verlag, Berlin, 2005. A global geometric and
probabilistic perspective, Mathematical Physics, III.\
\bib{1BGW}
C. Bonatti, S. Gan and L. Wen, On the existence of non-trivial
homoclinic class, {\it preprint} (2005)\
\bib{1BGV}
C. Bonatti, N. Gourmelon and T. Vivier, Perturbation of the
derivative along periodic orbits, {\it preprint} (2004)\
\bib{1Co}
C. Conley, Isolated invariant sets and Morse index, CBM Regional
Conference Series in Mathematics, {\bf 38}, AMS Providence,
R.I.,(1978).\
\bib{1C1}
S. Crovisier, Periodic orbits and chain transitive sets of $C^1$
diffeomorphisms, {\it preprint} (2004).\
\bib{1C2}
S. Crovisier, Birth of homoclinic intersections: a model for the
central dynamics of partial hyperbolic systems, {\it
preprint} (2006).\
\bib{1F}
J. Franks, Necessary conditions for stability of diffeomorphisms,
{\it Trans. Amer. Math. Soc.}, {\bf 158} (1977), 301-308.\
\bib{1G1}
S. Gan, Private talk.\
\bib{1G2}
S. Gan, Another proof for $C^1$ stability conjecture for flows,
{\it SCIENCE IN CHINA (Series A)} {\bf 41 No. 10} (October 1998)
1076-1082 \
\bib{1G3}
S. Gan, The Star Systems $X^*$ and a Proof of the $C^1
Omega$-stability Conjecture for Flows, {\it Journal of
Differential Equations}, {\bf 163} (2000) 1--17\
\bib{1GW1}
S. Gan and L. Wen, Heteroclinic cycles and homoclinic closures for
generic diffeomorphisms, {\it Journal of Dynamics and Differential
Equations,}, {\bf 15} (2003), 451-471.\
\bib{1GW2}
S. Gan and L. Wen, Nonsingular star flow satisfy Axion M and the
nocycle condition, {\it Ivent. Math.,} {\bf 164} (2006), 279-315.\
\bib{1H}
S. Hayashi, Connecting invariant manifolds and the solution of the
$C^1$ stability and $\Omega$-stable conjecture for flows, {\it
Ann. math.}, {\bf 145} (1997), 81-137.\
\bib{1HPS}
M. Hirsch, C. Pugh, and M. Shub, Invariant manifolds, volume 583
of {\it Lect. Notes in Math}. Springer Verlag, New york, 1977\
\bib{1L1}
Shantao Liao, Obstruction sets I, {\it Acta Math. Sinica}, {\bf 23} (1980), 411-
453.\
\bib{1L2}
Shantao Liao, Obstruction sets II, {\it Acta Sci. Natur. Univ. Pekinensis},
{\bf 2} (1981), 1-36.\
\bib{1L3}
Shantao Liao, On the stability conjecture, {\it Chinese Annals of Math.},
{\bf 1} (1980), 9-30.(in English)\
\bib{1L4}
Shantao Liao, An existence theorem for periodic orbits, {\it Acta Sci.
Natur. Univ. Pekinensis}, {\bf 1} (1979), 1-20.\
\bib{1L5}
Shantao Liao, Qualitative Theory of Differentiable Dynamical Systems,
{\it China Science Press}, (1996).(in English)\
\bib{1M1}
R. Ma\~{n}\'{e}, Quasi-Anosov diffeomorphisms and hyperbolic manifolds,
{\it Trans. Amer. Math. Soc.}, {\bf 229} (1977), 351-370.\
\bib{1M2}
R. Ma\~{n}\'{e}, Contributions to the stability conjecture, {\it Topology}, {\bf 17} (1978),
383-396.\
\bib{1M3}
R. Ma\~{n}\'{e}, An ergodic closing lemma, {\it Ann. Math.}, {\bf 116} (1982), 503-540.\
\bib{1M4}
R. Ma\~{n}\'{e}, A proof of the $C^1$ stability conjecture. {\it
Inst. Hautes Etudes Sci. Publ. Math.} {\bf 66} (1988), 161-210.\
\bib{1N1}
S. Newhouse, Non-density of Axiom A(a) on S 2. {\it Proc. A. M. S.
Symp pure math}, {\bf 14} (1970), 191-202, 335-347.\
\bib{1N2}
S. Newhouse, Diffeomorphisms with infinitely many sinks, {\it Topology}, {\bf 13}, 9-18, (1974).\
\bib{1PV}
J. Palis and M. Viana, High dimension diffeomorphisms displaying infinitely sinks, {\it Ann. Math.}, {\bf 140} (1994), 1-71.\
\bib{1Pl}
V. Pliss, On a conjecture due to Smale, {\it Diff. Uravnenija.}, {\bf 8} (1972),
268-282.\
\bib{1P}
C. Pugh, The closing lemma, {\it Amer. J. Math.}, {\bf 89} (1967), 956-1009.\
\bib{1PS1}
E. Pujals and M. Sambarino, Homoclinic tangencies and hyperbolicity
for surface diffeomorphisms, {\it Ann. Math.}, {\bf 151} (2000), 961-1023.\
\bib{1PS2}
E. Pujals and M. Sambarino, Density of hyperbolicity and
tangencies in sectional dissipative regions, {\it preprint}
(2005).\
\bib{1S}
J. Selgrade, Isolated invariant sets for flows on vector bundles, {\it Trans.
Amer. Math. Soc.}, {\bf 203} (1975), 259-390.\
\bib{1SS}
R. Sacker and G. Sell, Existence of dichotomies and invariant splittings
for linear differential systems, {\it J. Diff. Eq.}, {\bf 22} (1976), 478-496.\
\bib{1Sh}
M. Shub, Topological transitive diffeomorphisms in $T^4$, Lecture Notes
in Math. Vol. {\bf 206}, Springer Verlag, 1971.\
\bib{1W0}
L. Wen, On the $C^1$-stability conjecture for flows. {\it Journal
of Differential Equations,} {\bf 129}(1995) 334-357.\
\bib{1W1}
L. Wen, Homoclinic tangencies and dominated splittings, {\it Nonlinearity},
{\bf 15} (2002), 1445-1469.\
\bib{1W2}
L.Wen, A uniform $C^1$ connecting lemma, {\it Discrete and continuous dynamical systems}, {\bf 8} (2002), 257-265.\
\bib{1W3}
L. Wen, Generic diffeomorphisms away from homoclinic tangencies and
heterodimensional cycles, {\it Bull. Braz. Math. Soc. (N.S.)}, {\bf 35} (2004), 419-
452.\
\bib{1W4}
L. Wen, Selection of Quasi-hyperbolic strings, {\it preprint} (2006).\
\bib{1WG}
L. Wen and S. Gan, Obstruction sets, Obstruction sets, quasihyperbolicity
and linear transversality.(Chinese) {\it Beijing Daxue Xuebao
Ziran Kexue Ban}, {\bf 42} (2006), 1-10.\
\bib{1XW}
Z. Xia and L. Wen, $C^1$ connecting lemmas, {\it Trans. Amer.
Math. Soc.} {\bf 352} (2000), 5213-5230.\
\bib{1YGW}
D. Yang, S. Gan, L, Wen, Minimal Non-hyperbolicity and
Index-Completeness, {\it preprine} (2007).\
\bib{1Y1}
J. Yang, Ergodic measure far away from tangency, {\it
preprint} (2007).\
\bib{1Y2}
J. Yang, Lyapunov stable chain recurrent class, {\it preprint}
(2007).\
\bib{1Y4}
J. Yang, Aperiodic class, {\it preprint} (2007).\
\bib{1ZG}
Y. Zhang and S. Gan, On Ma\~{n}\'{e}'s Proof of the $C^1$
Stability Conjecture, {\it Acta mathematica Sinica, English
Series} Vol. {\bf 21}, No. {\bf 3}, June, 2005, 533-540.\
 \et

\end{document}